\documentclass{article}
\usepackage[utf8]{inputenc}
\usepackage{amsfonts,amsmath,amsthm,latexsym,amssymb,marvosym,esint,xfrac,bbm,amssymb}
\usepackage{graphics, epsfig}
\usepackage{color}
\usepackage{tikz}
\tikzstyle{mybox} = [draw=black, very thick, rectangle, rounded corners, inner ysep=5pt, inner xsep=5pt]

\usepackage[shortlabels]{enumitem}
\usepackage{appendix}

\usepackage{verbatim}
\usepackage{ulem}
\usepackage[makeroom]{cancel}
\usepackage[margin=1in]{geometry}
 \usepackage[usenames,dvipsnames]{pstricks}
 \usepackage{pst-grad} 
 \usepackage{pst-plot} 
\allowdisplaybreaks

 \numberwithin{equation}{section}
 \newtheorem{theorem}{Theorem}[section]
 
 \newtheorem{lemma}[theorem]{Lemma}
\newtheorem{proposition}[theorem]{Proposition}
 \newtheorem{remark}[theorem]{Remark}
 
 \newtheorem{corollary}[theorem]{Corollary}
 \newtheorem{definition}[theorem]{Definition}
 \theoremstyle{definition}

\def\R{\mathbb{R}}

\def\Q{\mathcal{Q}}

\def\K{\mathcal{K}}
\def\N{\mathbb{N}}

\def\I{\mathcal{I}}

\def\T{\mathcal{T}}
\def\S{\mathcal{S}}

\title{Harnack-type estimates and extinction in finite time for a class of anisotropic porous medium type equations\\
{\small  }}
\author{ Eurica Henriques {\&} Simone Ciani }

\begin{document}

\maketitle

\section{Abstract} 
In this work we are interested in the study of a class of anisotropic porous medium-type equations whose prototype is 
\begin{equation} \label{APME}
u_t =\sum_{i=1}^N \left( m_i u^{m_i-1} u_{x_i} \right)_{x_i}  \ , \qquad 0<m_1 \leq \cdots \leq m_N <1 \ ,
\end{equation}  
for which we derive several estimates, namely two Harnack-type inequalities; and, when considering the associated Dirichlet problem, we determine the finite time of extinction and thereby present a decay rate of extinction.

\noindent MSC-Mathematical Subject Classification- 
35K55, 35K67, 35B65, 35Q35, 35D30 \\

\noindent Key Words: Anisotropic Porous Medium -type Equations, Finite Time of Extinction, Rate of Extinction, Harnack-type Estimates

\section{Introduction}

\noindent Several physical phenomena, such as groundwater infiltration or heat radiation in plasmas, are described by nonlinear evolutionary differential equations of the form
\begin{equation} \label{PME}
u_t =\sum_{i=1}^N \left( m u^{m-1} u_{x_i} \right)_{x_i} =  \Delta \left( u^{m} \right) \ .
\end{equation}
These equations are the nonlinear version of the heat equation ($m=1$) and are known in the literature as the porous medium equation, when $m>1$; when $0<m<1$, the equation is called the fast diffusion equation (some authors, call it porous medium type equation). The mathematical challenges within their structure together with their physical relevance (due to the variety of applications) rouse the interest of several authors since the middle of the 20th century (see for instances \cite{Aron85}-\cite{DGVmono}-\cite{Vaz06}-\cite{Vaz12} and the references therein). It is known that the properties exhibit by the solution of \eqref{PME}  for $m\neq 1 $ do not hold when $m=1$: for $m>1$, the disturbances from the data propagate with finite speed, and, for $0<m<1$, solutions become zero (extinct) in finite time. These two properties were first obtained in \cite{BH80} (we refer to \cite{Sabinina62}, \cite{ZK50}, in the case of the p-Laplacian). Although the subject is not new, it is still deserving the attention of a wide community of mathematicians as can be seen in \cite{LWW07}, \cite{BGV12} and \cite{Porzio23}; and also in the recent works \cite{DFZ20}, \cite{BF21}, \cite{CMcCS23} and \cite{MNS23} where a doubly nonlinear parabolic equation is considered. All these works, although presenting different settings (for instance, the way the solution is defined), methods and approaches, rely closely on the isotropic character of the diffusion process. When in presence of electrorheological fluids, for instance when considering the flow of a barotropic gas in a nonhomogeneous anisotropic porous medium, the exponent associated with the diffusion is now a function $\gamma(x,t)$ of the thermodynamics, being the process modelled by 
\begin{equation} \label{electrorh}
u_t -\mathrm{div} \left(|u|^{\gamma(x,t)} Du \right)= f  \ . 
\end{equation}

\noindent Results on the existence, uniqueness and localization properties of the solutions to \eqref{electrorh} were obtained in \cite{AS05}; in \cite{HU06}-\cite{EH08}-\cite{EH21} a local regularity theory was developed (in \cite{AS15} and the references therein one can find a more complete description of the subject). 

\vspace{.2cm}

\noindent Consider now the water motion in an anisotropic porous medium. This particular feature changes the scenario: we are in the presence of an anisotropic phenomenon modelled by \eqref{APME}. In \cite{Songexist01}-\cite{Songuniq01}, the author proves the existence and uniqueness, respectively, of generalized solutions (also continuous by definition) to an initial-boundary value problem associated to an anisotropic porous medium equation possessing also singular anisotropic advection and strong absorption terms. In these works, the author stresses out the difficulties rising from the anisotropies of the diffusions as well as the singularities of the advections. In \cite{SongJian05}, the authors proved the existence of fundamental solutions (continuous by definition) to \eqref{APME}, for $0<m_i<1$, pointing out the need to take into account scaling techniques in each space direction and the construction of suitable supersolutions. In \cite{SongJian06}, the same authors establish the existence of an unique solution to a Cauchy problem with integrable initial data  as well as a comparison principle; the adopted setting considers the arithmetic mean $m={\sum_{i}^N m_i}/{N}$ to be in the supercritical range $m>(N-2)/N$, $m_N < (2+Nm)/N$ and $L^1$-regular solutions. Results on the local regularity of the weak solutions to \eqref{APME}, based on intrinsic isotropic scaling, were derived in \cite{EH11}, separately, for the degenerate case, $m_i>1$, and for the singular one, $0<m_i<1$. In the recent work \cite{FVV23}, the authors contribute with the analysis of self-similarity: they prove the existence of an unique fundamental solution of self-similar type and present, in terms of the family of self-similar fundamental solutions, the asymptotic behaviour of all finite mass solutions. Their work is developed for the prototype anisotropic fast diffusion equation \eqref{APME} defined along the stripe $\R^N \times \R^+$. 

\vspace{.2cm}

\noindent To the best of authors' knowledge and within the context presented here and explicited below, no results related to Harnack-type inequalities and extinction profile were obtained so far. This is precisely the scope of this work: not only to bring to light Harnack-type inequalities and extinction profile towards (singular) anisotropic porous medium type equations, but also (and no less significant) to do it in an wider context  where lower order terms are considered, so that absorption and advection may be taken into account. In fact, we will consider the following class of anisotropic evolution equations
\begin{equation} \label{gAPME} u_t -\mathrm{div} A(x,t, u, Du) = B(x,t,u,Du) \ , \qquad \mathrm{in} \quad  \Omega_T= \Omega \times (0,T]
\end{equation}
for positive $T$ and a bounded set $\Omega$ of $\R^N$, $N>2$, being $A= (A_1, \cdots, A_N)$ where $A_1, \cdots, A_N, B: \Omega_T 
\times\R^{N+1} \rightarrow \R$ are measurable functions satisfying the structure conditions 
\begin{equation} \label{SC1}
   A(x,t, u, Du) \cdot Du \geq C_o \sum_{i=1}^N m_i u^{m_i-1} |u_{x_i}|^2 -C^2 \sum_{i=1}^N u^{m_i+1} 
\end{equation}

\begin{equation} \label{SC2}
   |A(x,t, u, Du)|\leq C_1 \sum_{i=1}^N m_i u^{m_i-1} |u_{x_i}| +C \sum_{i=1}^N u^{m_i} 
\end{equation}

\begin{equation} \label{SC3}
  |B(x,t,u,Du)| \leq C \sum_{i=1}^N m_i u^{m_i-1} |u_{x_i}| + C^2 \sum_{i=1}^N u^{m_i} 
\end{equation}
for given positive constants $C_o$ and $C_1$, $C$ a nonnegative constant and $0<m_1 \leq \cdots \leq m_N <1$.

\subsection{Setting the framework}

\noindent To clarify notation, for a real valued function $u(x,t)$, being $(x,t) \in \Omega_T $ with $x=(x_1, \cdots, x_N)$, we consider the time and space derivatives as 

\[u_t= \frac{\partial u}{\partial t} \ , \qquad u_{x_i}= \frac{\partial u}{\partial x_i} \ , \ i=1, \cdots, N \ ,  \qquad Du= \left(u_{x_1}, \cdots, u_{x_N}\right) \ . \]

\vspace{.2cm} 

\noindent Within the setting \eqref{gAPME}-\eqref{SC1}-\eqref{SC2}-\eqref{SC3}, for $0<m_1 \leq \cdots \leq m_N <1$, we are in presence of a singular anisotropic differential equation: on the one hand, in every single direction, the modullus of ellipticity $u^{m_i-1}$ becomes unbounded at the points where $u$ vanishes - singular character -; on the other hand, the diffusion occurs differently in each space direction since there is an $i$-dependence on the exponents $m_i$ - anisotropic character. 

\vskip0.2cm 

\noindent This equation gives rise to two particular cases: when $C=0$, no lower order term is considered and we say the equation related to \eqref{gAPME} is homogeneous; if furthermore $C_o=C_1=1$, we recover the prototype anisotropic porous medium type equation \eqref{APME}.

\vskip0.2cm 

\noindent As mentioned before, the definition of solution taken at hands plays an important role when deriving estimates, properties and other results. In what follows we present what we mean by a local weak solution to \eqref{gAPME}-\eqref{SC1}-\eqref{SC2}-\eqref{SC3}.

\begin{definition} \label{defweaksol} A measurable function $u: \Omega_T \rightarrow \R$ is called a local weak sub(super)solution to \eqref{gAPME}-\eqref{SC1}-\eqref{SC2}-\eqref{SC3} if
\[ 
u \in C\left(0,T; L^2_{loc}(\Omega)\right) \quad \mathrm{and} \quad u^{\frac{m_i-1}{2}} u_{x_i} \in L^2_{loc}\left(0,T; L^2_{loc}(\Omega)\right),\quad  \forall i \in \{1,\dots, N\},
\] 
and, for all $K\subset \subset \Omega$ and for all $[t_1,t_2] \subset (0,T]$, the inequality
\begin{equation} \label{weaksol}
  \int_{K} u \varphi (x,t_2) \ dx - \int_{K} u \varphi (x,t_1) \ dx - \iint_Q u \varphi_t \ dx dt +  \iint_Q A(x,t, u, Du) \cdot D\varphi \ dx dt \leq (\geq) \iint_Q B(x,t,u,Du) \varphi \ dx dt
\end{equation}
holds true for all nonnegative test functions $ \varphi \in W^{1,2}_{loc}\left(0,T; L^2(K)\right) \cap L^2_{loc}\left(0,T; W^{1,2}_o(K)\right)$, being $Q=K \times [t_1, t_2]$.

\vspace{.2cm}

\noindent A function $u$ that is both a local weak sub and super solution is called a local weak solution to \eqref{gAPME}-\eqref{SC1}-\eqref{SC2}-\eqref{SC3}.

\end{definition}

\noindent Along the text, we will take test functions that depend on the solution itself, and thereby have (if any) low regularity in time (in general the time derivative of $u$ only makes sense in a distributional context). To present accurate proofs one should consider either the Steklov average of $u$, as proposed in \cite{DGVmono}, or the regularization $u^\star$ proposed in \cite{KL08} for the porous medium equation. Loosely speaking, this procedure allows to work freely with time derivatives of the regularized solution and then, due the regularity assumptions on $u$, pass to the limit. The procedure and limit process are very similar to the one already known for the isotropic equations. Hence, for the sake of simplicity and to keep the focus on the anisotropy, we decided to proceed in a formal fashion. For an accurate approach we refer to \cite{FHV21}, with the obvious changes.\newline

\subsection{Novelties and main results}

\noindent The study of the anisotropic porous medium type equation within this wider context is new and adds up several other difficulties to the anisotropy. Indeed, not only one has to deal with the anisotropic character of the differential equation, demanding and difficult by itself, but also one has to cope with the lower order terms: roughly speaking and just to illustrate what one has to deal with, when considering the prototype equation \eqref{APME} one can skip to take space derivatives of $u^{m_i}$ and just work with 
\[ \int_{K} u \varphi (x,t_2) \ dx - \int_{K} u \varphi (x,t_1) \ dx - \iint_Q u \varphi_t \ dx dt + \sum_{i=1}^N \iint_Q  u^{m_i} \varphi_{x_i x_i} \ dx dt =0 \ ,  \]
which is no longer possible in this large setting.  
\vskip0.1cm 

\noindent Unlikely to what happens to the porous medium type equation, where the transformation $v=u^m$ can be considered, in the anisotropic setting that is no longer possible; so a different nonlinear strategy has to be taken into account whose use might be adapted to a variety of other different situations. We would like also to mention a technical improvement that we think it may be useful in other contexts: the proof of the integral form of Harnack type inequality, at a certain stage, requires to work with some sort of energy estimates exhibiting negative powers. This turned out to be a  technical difficulty which we overcome by keeping the anisotropic estimate as it appeared and afterwards assume a certain condition on the supremum of $||u(\cdot, \tau)||_{L^1}$, for $0 \leq \tau \leq t$. For an insight on this simple but effective technique, we refer to the proof of Lemma \ref{lemma1} and later condition \eqref{condsup}, for the intrinsic anisotropic setting; for the standard setting, see Lemma \ref{lemma1isot} and \eqref{kondizionen}. 

\vskip0.2cm \noindent 
\noindent As said before, the anisotropy exhibited by the equation is quite a hard thing to deal with. How can one entangle this relevant and difficult feature along with the theory of Harnack-type inequalities? We realized that the adopted geometry plays a crucial role along the process (and there is a price to be paid), unlikely to what happens in the isotropic setting. Our local results will be presented and proved under two different geometries: for an interior point $(x_o,t_o) \in \Omega_T$, we consider 

\begin{itemize}
    \item[*] the standard one, meaning a geometry that keeps time and space independent, {\it{i.e.}} 
    \[Q(x_o,t_o)= \prod \left\{|x_i-x_{oi}| < \rho^{p/p_i} \right\} \times [t_o, t_o+ t] \  \]
More details about this geometry and results obtained under its choice can be found in Section \ref{S:isotgeom}.
    \item[*] an intrinsic anisotropic geometry, meaning a geometry for which time tangles within the cube's radius. More precisely, for fixed $\rho>0$ and $t>0$, we consider the intrinsic anisotropic cubes and cylinders
\begin{equation} \label{AnisoCubes}{\cal K}_{a\rho} (x_o) = \prod_{i=1}^N \left\{ |x_i-x_{oi}| < \left(\frac{t}{\rho^2}\right)^{\frac{m_i-m}{2(1-m)}} a \rho \right\}, \quad  {\cal Q}_{a\rho} (x_o,t_o) = {\cal K}_{a\rho} (x_o)\times [t_o,t_o+t] \ \ , \ \ a>0 \ .  \end{equation}
In Section \ref{S:anisotPME} this geometry and its choice will be made precise. 
\end{itemize}

\vspace{.2cm}

\noindent In what follows we present the main results stated within the anisotropic geometric framework (under a translation argument, we will consider $(x_o,t_o)=(0,0)$), starting with two Harnack-type inequalities and finalizing with results related to the existence of a finite time of extinction and its associated decay rate. Analogous results  will be derived (and proved) in Section \ref{S:isotgeom} where a standard isotropic geometry is undertaken.

\begin{theorem}{\bf [Integral form of a Harnack-type inequality]}
\label{integralharnack}
\vskip0.1cm 
\noindent Let $u$ be a nonnegative, local weak solution to \eqref{gAPME}-\eqref{SC1}-\eqref{SC2}-\eqref{SC3} in $\Omega_T$. There exists a positive constant $\gamma$, depending on $N, C_o, C_1, m_N$, such that, for all cylinders ${\cal Q}_{2\rho} \subset \Omega_T$, either
\begin{equation}\label{condCrhoT1}
C \rho \left(\frac{t}{\rho^2}\right)^{\frac{m_i-m}{2(1-m)}} >1, \quad \mathrm{for} \ \mathrm{some} \ i=1, \cdots,N \ ,
\end{equation}
or 
\begin{equation} \label{estL1L1}
\sup_{0\leq \tau \leq t}  \int_{{\cal K}_{\rho}} u(x, \tau) \ dx  \leq \gamma \ \left\{ \inf_{0\leq\tau \leq t} \int_{{\cal K}_{2\rho}} u(x, \tau) \ dx + \left(\frac{t}{\rho^\lambda} \right)^{\frac{1}{1-m}} \right\} \ , 
\end{equation}
where $\lambda= N(m-1)+2$.
\end{theorem}

\begin{theorem} {\bf [$L^1_{loc}-L^{\infty}_{loc}$ Harnack-type estimate]}
\label{L1Linftest}\vskip0.1cm 
\noindent Let $u$ be a nonnegative, locally bounded, local weak solution to \eqref{gAPME}-\eqref{SC1}-\eqref{SC2}-\eqref{SC3} in $\Omega_T$, for $\displaystyle{m >{N-2}/{N}}$. There exists a positive constant $\gamma$, depending on $N, C_o, C_1, m_i$, such that, for all cylinders ${\cal Q}_{2\rho} \subset \Omega_T$, either
\begin{equation}\label{condCrho}
C \rho \left(\frac{t}{\rho^2}\right)^{\frac{m_i-m}{2(1-m)}} >1, \quad \mathrm{for} \ \mathrm{some} \ i=1, \cdots,N, 
\end{equation}
or 
\begin{equation} \label{estL1Linfty}
\sup_{\K{_{\rho/2}} \times [t/2,t]}u \leq \gamma \ t^{-\frac{N}{\lambda} } \left(\inf_{ 0 \leq \tau \leq t} \int_{\K_{2\rho}} u (x, \tau) \ dx \right)^{2/\lambda} + \gamma \left(\frac{t}{\rho^2} \right)^{\frac{1}{1-m}}  \ .\ 
\end{equation}
\end{theorem}

\noindent \begin{remark} The alternative condition \eqref{condCrhoT1} presented in both Theorem \ref{integralharnack} and Theorem \ref{L1Linftest} is related to the existence of lower order terms: it means that some projection of the set $\K_{\rho}$ along a coordinate axis has to be bigger than $1/C$ and implies an estimate for the integral terms related to the lower order terms. When equation \eqref{gAPME} is homogeneous ($C=0$) this alternative argument is void, while when all $m_i \equiv m$, that is, when equation  \eqref{gAPME} becomes the fast diffusion equation, it reduces to the classical alternative $C\rho >1$. 

\noindent If one adopts a standard geometry, the second terms on the right-hand sides of \eqref{estL1L1} and \eqref{estL1Linfty} is replaced by a sum where each component will depend on $m_i$ instead of $m$; see Section \ref{S:isotgeom} for a complete understanding.
\end{remark}

\vspace{.2cm}

\noindent Another novelty concerns the attainment of extinction in finite time, regardless of the geometry undertaken since one works with a rectangular bounded domain $\Omega \subset \R^N$ containing the origin. In fact, when considering the Dirichlet problem 

\begin{equation} \label{DP}
\left\{
\begin{array} {ll}
u_t - \mathrm{div} A_o(x,t,u,Du) =0 & \ , \ (x,t) \in  \Omega \times \R^+\\[.8em]
u(x,t) =0 & \ , \ (x,t) \in \partial \Omega \times (0,+\infty)\\[.8em]
u(x,0)=u_o(x) \geq 0 & \ , \ x \in \Omega
\end{array} \right.
\end{equation}
where the operator $A_o$ satisfies the homogeneous structure conditions,
\begin{equation} \label{SC1H}
   A_o(x,t, u, Du) \cdot Du \geq C_o \sum_{i=1}^N m_i u^{m_i-1} |u_{x_i}|^2  \ , \quad C_o>0
\end{equation}

\begin{equation} \label{SC2H}
   |A_o(x,t, u, Du)|\leq C_1 \sum_{i=1}^N m_i u^{m_i-1} |u_{x_i}|  \ , \quad C_1>0 \ ,
\end{equation}
assuming that \eqref{DP} has a unique solution (see for the special cases \cite{AS05},\cite{AS15}, \cite{SongJian05},\cite{SongJian06}), this solution vanishes in a time $T^*$ that can be quantified in terms of some $L^r$-norm of the initial datum $u_o$.
\begin{theorem}{\bf [Finite time of extinction in bounded domains]}\label{extfinitetime}
\vskip0.1cm 
\noindent Let $\Omega$ be a rectangular bounded domain in $\R^N$, $N>2$. Assume  $u$ is the unique nonnegative, bounded weak solution to \eqref{DP}-\eqref{SC1H}-\eqref{SC2H}, where $u_o \in L^\infty(\Omega)$. Then there exits a finite positive time $T^\star$, depending on $C_o, m_1,m_N, m, N$, such that
\[u(x,t) \equiv 0 \ , \quad \mathrm{for} \ \ \mathrm{all} \ \ t \geq T^\star \ . \]
Moreover, letting $\beta_1= [N(m-1)+2 +2m_N]/[N(m_N+1)^2]$, the extinction time $T^*$ has an upper bound:
\begin{equation} \label{Tstar}
    T^\star \leq  \left\{
    \begin{array}{ll}
\gamma_1 \ ||u_o||_{L^{N(1-m)/2}(\Omega)}^{1-m} & \ ,  \ m_N <  \frac{N(1-m)}{2} -1 \quad \text{and} \quad 0<m<\frac{N-2}{N} \\[.8em]
 \gamma_2  \ | \Omega|^{- \beta_1} \ ||u_o||_{L^{m_N+1}(\Omega)}^{1-m} &
 \ , \ m_N \geq \frac{N(1-m)}{2} -1 \quad \text{and} \quad0<m<\frac{N-2}{N} \\[.8em]
\gamma_2  \ | \Omega|^{- \beta_1} \ ||u_o||_{L^{m_N+1}(\Omega)}^{1-m} & \ , \ \frac{N-2}{N} \leq m <1,
    \end{array}
    \right.
\end{equation} \noindent where the positive constants $\gamma_1$ and $\gamma_2$ only depend on $C_o$, $m_1$, $m_N$, $m$, $N$.
\end{theorem}
\noindent From this result one is able to derive the decay rate of extinction that relies on Harnack-type estimates and therefore on the adopted geometry (see Section\ref{S:isotgeom} for details on this topic when one chooses to work within the standard geometry). 

\begin{theorem}\label{decayinfty}
{\bf [Decay rate of extinction]}
\vskip0.1cm 
\noindent
In the setting of Theorem \ref{extfinitetime}, let $m>(N-2)/N$ be supercritical. Then, there exists a positive constant $\gamma$, depending upon $C_o,m_1, m_N, N,m$, such that
\[ ||u(\cdot, t)||_{L^{\infty}(\K_\rho)}  \leq \gamma \  \left(\frac{T^\star -t}{\rho^2} \right)^{\frac{1}{1-m}} \ \ , \qquad \mathrm{for} \ \ \mathrm{all} \ \  {T^\star}/{2} < t <T^\star , \]
being $T^\star$ be a finite time of extinction and assuming $\K_{\rho}= \prod_{i=1}^N \left\{|x_i|< \left(\frac{2(T^*-t)}{\rho^2}\right)^{(m_i-m)/[2(1-m)]} \rho\right\}\subset \Omega$.
\end{theorem}

\subsection{Adopted text structure}

\noindent In Section \ref{S:knownresults}, to keep the text as self-contained as possible, we present several known results; in Section \ref{S:anisotPME}, we present an intrinsic anisotropic geometry within it \eqref{gAPME} behaves as the (isotropic) porous medium type equation with respect to the arithmetic average $m=\sum_{i=1}^N m_i /N $. In Section \ref{S:integralHarnack} and Section \ref{S:l1linftHarn}, working on the framework of intrinsic anisotropic geometry, we derive an integral Harnack-type inequality and a $L^1$-$L^\infty$ Harnack-type inequality, respectively. Still under the framework of intrinsic anisotropic geometry, Section \ref{S:Lrbackaniso} comprehends the presentation and proof of a local $L^r$ estimate backward in time, while Section \ref{S:extinction} is devoted to the study of the extinction profile. Finally, in Section \ref{S:isotgeom}, all the previous results are revisited and proved within the context of the standard geometry.


\section{Auxiliary results} \label{S:knownresults}

As it is well known, for the isotropic porous medium equation (see for instances \cite{DGVmono}), the proof of Harnack-type inequalities relies on two important estimates: the H\"{o}lder inequality and the Sobolev-Nirenberg embedding. In the present setting, that is, in the context of anisotropic diffusion equations, we will also need to have this kind of embeddings at hand. 

\vspace{.2cm}
\noindent While the first result, originally investigated by Troisi \cite{Troisi}, provides an anisotropic elliptic embedding, the second one (see \cite{DMV}), gives us the expected embedding in the parabolic framework. Although these two results are given in a wider fashion (see also \cite{Besov}), for any $1< p_i < \infty$  ($i=1, \cdots,N$),  for our purposes it suffices to present them in the case of $L^2$-norms.

\begin{proposition}{[Anisotropic Sobolev-Troisi embedding]} 
\label{embeddingellip} \vskip0.1cm \noindent 
Let $\Omega\subseteq \R^N$ be a rectangular domain, for $N>2$, and consider $\alpha_{i}>0$,  $i=1,\dots, N$. Define 
\[ \alpha = \sum_{i=1}^N \alpha_i \qquad \mathrm{and} \qquad 2^{\star}_{\alpha} = 2^\star \frac{\alpha}{N} ,  \]
where $2^\star= {2N}/{N-2}$ stands for the Sobolev exponent . Then, there exists a positive constant $C=C(N, \alpha)$ such that, for every $u\in W^{1,2}_o (\Omega)$,
\begin{equation}
\label{EE}
||u||_{L^{2^{\star}_{\alpha}}(\Omega)} \leq C \prod_{i=1}^N || (|u|^{\alpha_i})_{x_i}||_{L^2(\Omega)}^{ 1/\alpha}
\end{equation}
\end{proposition}

\vspace{.2cm}

\begin{proposition}{[Anisotropic Gagliardo-Sobolev-Nirenberg embedding]} 
\label{embeddingparab} \vskip0.1cm \noindent 
Let $\Omega\subseteq \R^N$ be a rectangular domain, $N>2$, $\alpha_{i}>0$, $i=1,\dots, N$, and $\sigma\in [1, 2_{\alpha}^{*}]$, being $2_{\alpha}^{*}$ as before. \newline For any number $\theta\in [0,\, (N-2)/{N}]$, define 
\[
q=q(\theta, \alpha)=\theta \,  2^{*}_{\alpha}+\sigma\, (1-\theta).
\]
Then there exists a positive constant $\tilde{c}=\tilde{c}(N, \alpha, \theta, \sigma)$ such that
\begin{equation}
\label{PS}
\iint_{\Omega_{T}}|u|^{q}\, dx\, dt\leq \tilde{c}\, T^{1-\theta\, \frac{N}{N-2}}\left(\sup_{t\in [0, T]}\int_{\Omega}|u|^{\sigma}(x, t)\, dx\right)^{1-\theta} \ \prod_{i=1}^{N}\left(\iint_{\Omega_{T}}| \left(|u|^{\alpha_{i}}\right)_{x_i}|^{2}\, dx\, dt\right)^{\frac{\theta}{N-2}},
\end{equation} 
for any $u\in L^{1}(0, T; W^{1,1}_{0}(\Omega))$ for which the right hand side is finite (otherwise \eqref{PS} is trivially true). 
\end{proposition}

\vspace{.5cm} 

\noindent The two nonlinear iteration Lemmata that we are about to describe can be found in \cite{DB} ($c>1$) or \cite{CVV} ($c>0$), and are related to some properties of sequences of numbers $(Y_n)_n$, $n\in \mathbb{N}_o$. The first one concerns the geometric convergence of $(Y_n)_n$, while the second one turns the qualitative information on an equibounded sequence $(Y_n)_n$ into qualitative information on its first element $Y_o$.

\begin{lemma}{[Fast geometric convergence Lemma]\qquad }
\label{fastgeomconv}\vskip0.1cm 
\noindent Let $(Y_n)_n$ be a sequence of positive numbers verifying
\[ Y_{n+1} \leq c b^n \ Y_n^{1+\alpha} \ , \]
being $c>0$, $b>1$ and $\alpha>0$ given numbers. If $Y_o \leq c^{-1/\alpha} \ b^{-1/\alpha^2}$ then $Y_n \rightarrow 0$, as $n \rightarrow \infty$.
\end{lemma}

\vspace{.2cm} 

\begin{lemma}{[Interpolation Lemma]}
\label{iteration}\vskip0.1cm \noindent 
If we have a sequence of equibounded positive numbers $\{Y_n\}$ such that
\begin{equation}\label{sqn}
    Y_n \leq \epsilon Y_{n+1}+ c b^n, \quad c,b>1, \quad \epsilon \in (0,1),
\end{equation} \noindent then there exists $\gamma = \gamma(b,c)>0$ such that 
\[Y_0 \leq \gamma \ c.\]
\end{lemma}


\section{ Intrinsic anisotropic geometry} \label{S:anisotPME}

\noindent In this Section we introduce the  intrinsic anisotropic geometry for which the results in the next three Sections are presented (and proved). For this purpose, let $\rho>0$ and $t>0$ be fixed and consider the arithmetic mean $m:= {\sum_{i=1}^N m_i}/{N}$. We consider the anisotropic cubes 
and the corresponding cylinders centered at the origin 
\begin{equation} \label{AnisoCubes}{\cal K}_{a\rho} = \prod_{i=1}^N \left\{ |x_i| < \left(\frac{t}{\rho^2}\right)^{\frac{m_i-m}{2(1-m)}} a \rho \right\} \qquad \mathrm{and} \qquad  {\cal Q}_{a\rho}  = {\cal K}_{a\rho} \times [0,t] \ \, \quad \mathrm{for} \ \  \ a>0 \ .  
\end{equation}
Observe that, although in each space direction we have intervals with different lengths (some are very small and others are big, depending on the ratio $(t/\rho^2)$ and on the difference $m_i-m$), this geometry {\it preserves the volume}:
\[
|{\cal K}_\rho|= 2^N \rho^N = |K_\rho|,
\]
being $K_\rho= \left\{ |x| < \rho\right \}$ the usual cube in $\R^N$ with edge $2\rho$.

\vspace{.3cm}

\noindent This type of anisotropic geometry will play a crucial role when deriving homogeneous estimates for $u$. If not for anything else, this justifies (in some sense) its choice and use. However one can go further on and understand it within a diverse but related subject: self-similar solutions and geometries for which the energy estimates (fundamental tools in the theory undertaken) are invariant. We refer to \cite{Barenblatt-scaling1}, \cite{CGV22}, \cite{CSV22} and \cite{FVV23} one some account on this topic.

\noindent To ease notation and have a glimpse on this topic, consider the prototype anisotropic porous medium type equation \eqref{APME}. Let $M,L_i,T \in \R^+ $ and consider the change of variables 
\[ y=(y_1, \cdots, y_N) = (L_1 x_1, \cdots, L_N x_N) \quad \mathrm{and } \quad z= Tt \]
and the function
\[v(x, t)= M^{-1}u(y,z) \ . \]
If we ask both $u,v$ to be solutions to \eqref{APME}, they are called self-similar, and this determines a relation between the coefficients $M, L_i, T$. In fact, observe that, since
\[ v_t(x,t)= M^{-1} T \ u_z(y, z) \qquad \mathrm{and} \qquad v_{x_i}(x,t)= M^{-1} L_i \  u_{y_i}(y,z),\]

\[ u_z = \sum_{i=1}^N \left( m_i u^{m_i-1} u_{y_i} \right)_{y_i}  \ \ \Longrightarrow \ \ v_t= \sum_{i=1}^N \left[ M^{m_i-1} L_i^{-2} T \ \left(m_i v^{m_i-1} v_{x_i} \right)_{x_i} \right]  \ . \]
The homogeneity is reestablished once we take, for all $i=1, \cdots, N$,
\begin{equation} \label{self-similar-rule}
M^{m_i-1}T=L_i^2 \quad \Rightarrow \quad L_i=M^{\frac{m_i-1}{2}}T^{\frac{1}{2}}.\end{equation}

\noindent Let us consider $T^{\frac{1}{2}}=M^{\frac{1-m}{2}} A$, for some positive real number $A$. Then \eqref{self-similar-rule} can be written as \[L_i= M^{\frac{m_i-m}{2}} A, \qquad \forall i =1,\dots, N.\]
Therefore from various possible choices, that usually take into account the conservation of mass (see \cite{CGV22} and \cite{FVV23}), we let $Q_1= K_1 \times [0,1]$ and focus on the self-similar transformations 
\begin{enumerate}
    \item $v_1(x,t)= M^{-1} u (  M^{\frac{m_1-1}{2}} T^{\frac{1}{2}}x_1, \cdots, M^{\frac{m_N-1}{2}} T^{\frac{1}{2}}x_N, T t)= M^{-1} u (\T_M^T(x,t))$, with
    \[ \T_M^T( Q_1)= \prod_{i=1}^N \bigg{\{}|x_i|< M^{\frac{m_i-1}{2}} T^{\frac{1}{2}} \bigg{\}} \times [0,T] ,\]

    \item $v_2(x,t)= M^{-1} u ( A M^{\frac{m_1-m}{2}}  x_1, \cdots, A M^{\frac{m_N-m}{2}}  x_N,  A^2 M^{1-m}t)= M^{-1} u (\T_M^A(x,t))$ with 
    \[\T_M^A( Q_1)= \prod_{i=1}^N \bigg{\{}|x_i|< A M^{\frac{m_i-m}{2}}  \bigg{\}} \times [0,A^2 M^{1-m}].\]
    
\end{enumerate}

\noindent In the geometries given by $\T_M^T$ and $\T_M^A$ the energy estimates are invariant (see \cite{CGV22}). In addition, there are two interesting features to be considered: in $\T_M^T$ and for small $M$, stretched intervals are considered in every single space direction (since all $ 0<m_i<1$) which is consistent with the singular character of the equation, but the volume is not preserved. As for $\T_M^A$, the space intervals do not exhibit the same stretched behavior, while nevertheless the volume is preserved: this last property is the best suited for our aim. This discussion motivates the following introduction of geometry:  for fixed $\rho, t >0$, take $A= \rho $ and $M= \left({t}/{\rho^2} \right)^{\frac{1}{1-m}}$, then we obtain the intrinsic cylinders \eqref{AnisoCubes} as

\begin{equation} \label{geom}
\T_M^A(Q_1 ) = \K_{\rho}\times [0,t]\ . \end{equation}

\vspace{.5cm}

\noindent In the three following Sections, we consider that $\rho$ and $t$ are fixed positive real numbers such that 
\[\Q_{4\rho}= {\cal K}_{4\rho} \times [0,t] \subset \Omega_T \ .\]



\section{Proving the integral form of a Harnack-type inequality} \label{S:integralHarnack}

In this Section we prove Theorem \ref{integralharnack} with the help of the following result.

\begin{lemma}\label{lemma1}
\vskip0.1cm 
\noindent Let $u$ be a nonnegative local weak supersolution to \eqref{gAPME}-\eqref{SC1}-\eqref{SC2}-\eqref{SC3} in $\Omega_T$. Let $\alpha \in \R$ be such that $0< \alpha +1 <1$, $0<\alpha +m_i <1$ and $0<m_i-\alpha <1$. There exists a positive constant $\gamma$, depending on $N, C_o, C_1,m_1, \alpha$, such that, for all cylinders ${\cal Q}_{(1+\sigma)\rho} \subset \Omega_T$, for all $\sigma \in (0,1)$, either
\begin{equation}\label{condCrhoL1}
C \rho \left(\frac{t}{\rho^2}\right)^{\frac{m_i-m}{2(1-m)}} >1, \quad \mathrm{for} \ \mathrm{some} \ i=1, \cdots, N   \ ,
\end{equation}
or 
\begin{equation} \label{estL1}
\sum_{i=1}^N \iint_{{\cal Q}_{\rho}} u^{m_i+\alpha-2} |u_{x_i}|^2  \ dx d\tau \leq \frac{\gamma}{\sigma^2} \ \left\{ 1+ \sum_{i=1}^N \left( \left(\frac{t}{\rho^2}\right)^{\frac{1}{1-m}} \frac{\rho^N}{S} \right)^{1-m_i} \right\} 
\S_\sigma^{1+\alpha} \ \rho^{-N\alpha} \ ,
\end{equation}
where $\displaystyle{S_\sigma = \sup_{0 \leq \tau \leq t } \int_{{\cal K}_{(1+\sigma)\rho}} u  \ dx }$ and $\displaystyle{S = \sup_{0 \leq \tau \leq t } \int_{{\cal K}_{\rho}} u  \ dx }$.
\end{lemma}

\begin{proof}
   In the weak formulation \eqref{weaksol}, consider $\varphi=u^\alpha \xi(x)$, being 
    \[\xi(x)=\prod_{i=1}^N \xi_i(x_i)\]
    a smooth cutoff function defined in ${\cal K}_{(1+\sigma)\rho}$ and verifying $\xi=1$ in ${\cal K}_{\rho}$, $\xi=0$ outside ${\cal K}_{(1+\sigma)\rho}$ such that
    \[|\xi_i^\prime| \leq \frac{1}{\sigma \rho} \left(\frac{t}{\rho^2}\right)^{\frac{m-m_i}{2(1-m)}}    
 \ , \quad i=1, \cdots, N \, \]
and integrate over $ {\cal Q}_{(1+\sigma)\rho}$. Recalling that $\alpha$ is negative, we then have 
\begin{eqnarray*}
0 & \leq & \int_{{\cal K}_{(1+\sigma)\rho}} u^{\alpha+1}(x,t) \ \xi(x) \ dx - \int_{{\cal K}_{(1+\sigma)\rho}} u^{\alpha+1}(x,0) \ \xi(x) \ dx - \iint_{{\cal Q}_{(1+\sigma)\rho}} u \left(u^{\alpha}\right)_t \ \xi \ dx d\tau\\
& & + \iint_{{\cal Q}_{(1+\sigma)\rho}} A \cdot  D\left(u^{\alpha} \ \xi \right) \ dx d\tau\\
& & -  \iint_{{\cal Q}_{(1+\sigma)\rho}} B \  u^\alpha \xi \ dx d\tau\\
& = & I_1 + I_2 + I_3.
\end{eqnarray*}
Observe that
\begin{eqnarray*}
I_1 & = & \frac{1}{\alpha+1} \int_{{\cal K}_{(1+\sigma)\rho}} u^{\alpha+1}(x,t) \xi(x) \ dx - \frac{1}{\alpha+1} \int_{{\cal K}_{(1+\sigma)\rho}} u^{\alpha+1}(x,0) \xi(x) \ dx \\
& \leq & \frac{2}{\alpha+1}  \sup_{0 \leq \tau \leq t} \left(\int_{{\cal K}_{(1+\sigma)\rho}} u^{\alpha+1}(x,\tau) \ dx \right)^{\alpha +1} |{\cal K}_{(1+\sigma)\rho}|^{1-(\alpha+1)}\\
& \leq &\frac{2^{1+2N|\alpha|}}{\alpha+1} S_\sigma^{\alpha +1} \ \rho^{-N \alpha} \ . 
\end{eqnarray*}
Now, by recalling the structure conditions \eqref{SC1}- \eqref{SC2}-\eqref{SC3}, applying (twice) Cauchy's inequality to each $i$th-term, noting that
\[\iint_{{\cal Q}_{(1+\sigma)\rho}} u^{m_i+\alpha} \ dx d\tau  \leq 
t \ S_\sigma ^{m_i + \alpha} \ (2\rho)^{N(1-(m_i+\alpha))}
\]
and that \eqref{condCrhoL1} is violated, we get
\begin{eqnarray*}
I_2 & = & \alpha \sum_{i=1}^N  \iint_{{\cal Q}_{(1+\sigma)\rho}} A_i u_{x_i}  \ u^{\alpha-1} \  \xi \ dx d\tau + \sum_{i=1}^N  \iint_{{\cal Q}_{(1+\sigma)\rho}} A_i \  u^\alpha \ \xi_{x_i} \ dx d\tau \\
& \leq & 2\alpha C_o \sum_{i=1}^N m_i  \iint_{{\cal Q}_{(1+\sigma)\rho}} u^{m_i+\alpha-2} |u_{x_i}|^2 \xi^2 \ dx d\tau \\
& & + \frac{1}{(\sigma \rho)^2}  \sum_{i=1}^N  \left\{ -\alpha (C\rho)^2 -  \frac{C_1^2}{ \alpha C_o}  \left(\frac{t}{\rho^2}\right)^{\frac{m-m_i}{1-m}} + C \rho \left(\frac{t}{\rho^2}\right)^{\frac{m-m_i}{2(1-m)}} \right\} \iint_{{\cal Q}_{(1+\sigma)\rho}} u^{m_i+\alpha} \ dx d\tau \\
& \leq & 2\alpha C_o \sum_{i=1}^N m_i  \iint_{{\cal Q}_{(1+\sigma)\rho}} u^{m_i+\alpha-2} |u_{x_i}|^2 \xi^2 \ dx d\tau \\
& & + \frac{\gamma(N)}{|\alpha| \sigma^2}  \sum_{i=1}^N  \left\{ (C \rho)^2  \left(\frac{t}{\rho^2}\right)^{\frac{m_i-m}{1-m}} + \frac{C_1^2}{ C_o}  + C \rho  \left(\frac{t}{\rho^2}\right)^{\frac{m_i-m}{2(1-m)}} \right\} \left( 
\left(\frac{t}{\rho^2}\right)^{\frac{1}{1-m}}
\frac{\rho^N}{S_\sigma} \right)^{1-m_i}  \ S_\sigma^{\alpha +1} \ \rho^{-N \alpha}   \\
& \leq & 2\alpha C_o \sum_{i=1}^N m_i  \iint_{{\cal Q}_{(1+\sigma)\rho}} u^{m_i+\alpha-2} |u_{x_i}|^2 \xi^2 \ dx d\tau \\
& & +\frac{\gamma(N, C_o,C_1)}{|\alpha| \sigma^2}  \sum_{i=1}^N  \left( 
\left(\frac{t}{\rho^2}\right)^{\frac{1}{1-m}}
\frac{\rho^N}{S_\sigma} \right)^{1-m_i}  \ S_\sigma^{\alpha +1} \ \rho^{-N \alpha} \ ; 
\end{eqnarray*}
and 
\begin{eqnarray*}
I_3 & \leq & - \alpha C_o \sum_{i=1}^N m_i \iint_{{\cal Q}_{(1+\sigma)\rho}} u^{m_i+\alpha-2} |u_{x_i}|^2 \xi^2 \ dx d\tau \\
& & + \frac{C^2}{|\alpha |\sigma^2} \frac{C_o+1}{C_o} \sum_{i=1}^N  \iint_{{\cal Q}_{(1+\sigma)\rho}} u^{m_i+\alpha} \ dx d\tau\\
& \leq & - \alpha C_o  \sum_{i=1}^N  m_i \iint_{{\cal Q}_{(1+\sigma)\rho}} u^{m_i+\alpha-2} |u_{x_i}|^2 \xi^2 \ dx d\tau \\
& & +  \frac{\gamma(N)}{\alpha \sigma^2} \frac{C_o+1}{C_o} \sum_{i=1}^N  \left( 
\left(\frac{t}{\rho^2}\right)^{\frac{1}{1-m}}
\frac{\rho^N}{S_\sigma} \right)^{1-m_i}  \ S_\sigma^{\alpha +1} \ \rho^{-N \alpha} \ .
\end{eqnarray*}
Gathering all these estimates, and noticing that $S \leq S_\sigma$, we arrive at
\[\sum_{i=1}^N  \iint_{{\cal Q}_{(1+\sigma)\rho}} u^{m_i+\alpha-2} |u_{x_i}|^2 \xi^2 \ dx d\tau \leq \frac{\gamma(N,C_o, C_1, m_1)}{|\alpha|^2 (\alpha +1)} \frac{1}{\sigma^2} \left\{ 1+ \sum_{i=1}^N  \left( 
\left(\frac{t}{\rho^2}\right)^{\frac{1}{1-m}}
\frac{\rho^N}{S} \right)^{1-m_i}  \right\} \ S_\sigma^{\alpha +1} \ \rho^{-N \alpha} \ 
\]
which completes the proof.
\end{proof}


\subsection*{Proving Theorem \ref{integralharnack}}

\noindent Let $\rho>0$ and $t>0$ be fixed and construct the increasing sequence of anisotropic cubes
\[\K_{n}= \prod_{i=1}^N \left\{|x_i|< \bigg(\frac{t}{\rho^2} \bigg)^{\frac{m_i-m}{2(1-m)}} \rho_n \right\}\]
where 
\[\rho \leq \rho_n= \rho  \sum_{j=0}^n 2^{-j}  < 2 \rho \quad \mathrm{and} \quad \rho_{n+1}=(1+\sigma_n) \rho_n  \ \Rightarrow \sigma_n \geq \frac{1}{2^{n+1}} \ .\]
In the weak formulation \eqref{defweaksol} take $\varphi= \xi(x)$, being $\displaystyle{\xi(x)= \prod_{i=1}^N \xi_i(x_i)}$ a smooth cutoff function that: equals $1$ in $\K_{n}$, vanishes outside $\K_{n+1}$, therefore satisfying
\[|\xi_i^{\prime} | \leq \bigg(\frac{t}{\rho^2} \bigg)^{\frac{m-m_i}{2(1-m)}} \frac{2^{n+1}}{\rho} \ , \]
and consider the integration of \eqref{defweaksol} over $\K_{n+1} \times [\tau_1, \tau_2] \subset \K_{2 \rho} \times [0,t]$. We then get, 
\begin{eqnarray} \label{firstestintharnack}
\int_{\K_{n}}  u(x, \tau_1) \, dx & \leq  &\int_{\K_{n+1}} u(x, \tau_2)  \, dx \nonumber \\ & & +C_1 \frac{2^{n+1}}{\rho}  \sum_{i=1}^N m_i \bigg(\frac{t}{\rho^2} \bigg)^{\frac{m-m_i}{2(1-m)}} \int_{\tau_1}^{\tau_2}\int_{\K_{n+1}}  u^{m_i-1} |u_{x_i} | \, dxd\tau \nonumber \\
& & + C \sum_{i=1}^N m_i \int_{\tau_1}^{\tau_2}\int_{\K_{n+1}}u^{m_i-1} |u_{x_i} |\, dxd\tau\nonumber \\
& & + C \frac{2^{n+1}}{\rho}  \sum_{i=1}^N \bigg(\frac{t}{\rho^2} \bigg)^{\frac{m-m_i}{2(1-m)}} 
\int_{\tau_1}^{\tau_2}\int_{\K_{n+1}}  u^{m_i}  \, dxd\tau \nonumber \\
& & + C^2 \sum_{i=1}^N \int_{\tau_1}^{\tau_2}\int_{\K_{n+1}} u^{m_i} \, dxd\tau  \ .
\end{eqnarray} 
Now we choose $\tau_2 \in [0,t]$ such that
\[ \int_{\K_{2\rho}} u (x, \tau_2) \ d\tau = 
\inf_{0\leq \tau \leq t}  \int_{\K_{2\rho}} u(x,\tau) \ dx = \I \]
 and denote  $\displaystyle{S_n=\sup_{0\leq \tau \leq t}  \int_{\K_{n}} u(x, \tau) \ dx }$ and $\Q_{n}= \K_n \times [0,t]$. Thereby the previous inequality \eqref{firstestintharnack} now reads
 \begin{eqnarray*}
     S_n & \leq & \I + C_1 \frac{2^{n+1}}{\rho}  \sum_{i=1}^N  \bigg(\frac{t}{\rho^2} \bigg)^{\frac{m-m_i}{2(1-m)}} \iint_{\Q_{n+1}}  u^{m_i-1} |u_{x_i} | \, dxd\tau  + C \sum_{i=1}^N  \iint_{\Q_{n+1}} u^{m_i-1} |u_{x_i} |\, dxd\tau \\
& & + C \frac{2^{n+1}}{\rho}  \sum_{i=1}^N \bigg(\frac{t}{\rho^2} \bigg)^{\frac{m-m_i}{2(1-m)}} 
\iint_{\Q_{n+1}}u^{m_i}  \, dxd\tau 
 + C^2 \sum_{i=1}^N \iint_{\Q_{n+1}} u^{m_i}  \, dxd\tau  \ .
 \end{eqnarray*}
The terms evolving the directional space derivatives of $u$ are bounded from above with
\begin{eqnarray*}
\iint_{\Q_{n+1}}  u^{m_i-1} |u_{x_i} | \, dxd\tau  & = &
\iint_{\Q_{n+1}} \left( u^{\frac{m_i-1+\alpha -1}{2}} |u_{x_i} | \right) \left( u^{\frac{m_i-\alpha}{2}} \right) \, dxd\tau \\
& \leq & \left( \iint_{\Q_{n+1}} u^{m_i+\alpha -2} |u_{x_i} |^2 \, dxd\tau \right)^{1/2} \left( \iint_{\Q_{n+1}} u^{m_i-\alpha}  \, dxd\tau \right)^{1/2} \\
& \leq & \gamma \left(2^n  \ S_{n+1}^{(1+\alpha)/2} \ \rho^{-N\alpha/2} \right) \  \ \left(S_{n+1}^{(m_i-\alpha)/2} \ \rho^{N(1-m_i+\alpha)/2}  \  t^{1/2} \right)\\
& = & \gamma \ 2^n \ S_{n+1}^{(m_i+1)/2} \ \rho^{N(1-m_i)/2} \ t^{1/2}  \ .
\end{eqnarray*}
These estimates were obtained by: means of H\"{o}lder's inequality;  applying Lemma \ref{lemma1} to the pair of cylinders $\Q_n \subset \Q_{n+1}$;  recalling that \eqref{condCrhoT1} is not valid;  assuming, without loss of generality, that
\begin{equation} \label{condsup}
S= S_{o} > \left( \frac{t}{\rho^\lambda}\right)^{\frac{1}{1-m}} \ , 
\end{equation}
(otherwise, the results come immediately) and noting that, since $0< m_i-\alpha <1$,
\[\iint_{\Q_{n+1}}u^{{m_i-\alpha}} \, dxd\tau \leq \gamma(N) \ S_{n+1}^{m_i-\alpha} \ \rho^{N(1-m_i+\alpha)} \ t  \ .\]
Combining the estimate
\[\iint_{\Q_{n+1}}u^{m_i} \, dxd\tau \leq \gamma(N)  \ S_{n+1}^{m_i}\ \rho^{N(1-m_i)} \ t  \ ,\]
with the previous ones, while assuming \eqref{condCrhoT1} is violated, we get
\begin{eqnarray*}
S_n & \leq & \I + \gamma \ 2^{2n}  \sum_{i=1}^N \left\{ C_1 \bigg(\frac{t}{\rho^2} \bigg)^{\frac{m-m_i}{2(1-m)}} \left(\frac{ t}{\rho^2} \ \rho^{N(1-m_i)} \right)^{1/2} + C \rho \left(\frac{ t}{\rho^2} \ \rho^{N(1-m_i)} \right)^{1/2} \right\} S_{n+1}^{(m_i+1)/2}\\
& & + \gamma(N) \ 2^n \sum_{i=1}^N \left\{ C \rho \bigg(\frac{t}{\rho^2} \bigg)^{\frac{m-m_i}{2(1-m)}} \left(\frac{ t}{\rho^2}\right) \ \rho^{N(1-m_i)}  + (C \rho)^2 \left(\frac{ t}{\rho^2}\right) \ \rho^{N(1-m_i)} \right\} S_{n+1}^{m_i}\\
& \leq & \I + \gamma \ 2^{2n} \sum_{i=1}^N \left( \bigg(\frac{t}{\rho^2} \bigg)^{\frac{1-m_i}{1-m}}  \rho^{N(1-m_i)} \right)^{1/2}  S_{n+1}^{(m_i+1)/2}\\
& & + \gamma  \ 2^{2n}  \sum_{i=1}^N  \bigg(\frac{t}{\rho^2} \bigg)^{\frac{1-m_i}{1-m}}  \rho^{N(1-m_i)} S_{n+1}^{m_i} \\
& \leq &  \epsilon \ S_{n+1} + b^n \ \gamma(N,C_o,C_1,m_N,\epsilon) \left[ \left( \frac{t}{\rho^{\lambda}} \right)^{ \frac{1}{1-m}} + \I \right]\ , \quad b=2^{2/(1-m_N)}>1 \ .
\end{eqnarray*}
The last inequality was obtained by applying (at both summation terms) Young's inequality with $\epsilon_i=\epsilon$ in each $i$th-term. At this moment, we are two steps away to conclude the proof: firstly we iterate to get
\[S_o \leq \epsilon^n S_n + \sum_{k=0}^{n-1} (\epsilon b )^k \gamma  \left[ \left( \frac{t}{\rho^{\lambda}} \right)^{ \frac{1}{1-m}} + \I \right] \ ,
\]
and secondly we choose $\epsilon \in (0,1)$ such that $\epsilon b =1/2$. Since $\left(S_n\right)_n$ is equibounded, the result follows by letting $n \rightarrow \infty$.



\section{Proving local $L^1$-$L^{\infty}$ Harnack-type estimates} \label{S:l1linftHarn}

\noindent The main goal of this Section is to prove Theorem \ref{L1Linftest}, local $L^1$-$L^\infty$ Harnack-type estimates, for which one needs to derive local $L^r$-$L^\infty$ estimates, for $r \geq 1$. Namely,

\begin{proposition}{\bf [$L^r_{loc}-L^{\infty}_{loc}$ estimates]}
\label{LrLinftest}\vskip0.1cm 
\noindent Let $u$ be a nonnegative, locally bounded, local weak subsolution to \eqref{gAPME}-\eqref{SC1}-\eqref{SC2}-\eqref{SC3} in $\Omega_T$. Let $r \geq 1 $ be such that $\lambda_r= N(m-1)+2r >0$. There exists a positive constant $\gamma$, depending on $N, C_o, C_1, m_i$, such that, for all cylinders ${\cal Q}_{\rho} \subset \Omega_T$, either
\begin{equation}\label{condCrhoLrLinf}
C \rho \left(\frac{t}{\rho^2}\right)^{\frac{m_i-m}{2(1-m)}} >1, \quad \mathrm{for} \ \mathrm{some} \ i=1, \cdots,N
\end{equation}
or 
\begin{equation} \label{estLrLinfty}
\sup_{\K{_{\rho/2}} \times [t/2,t]}u \leq \gamma \ t^{-\frac{N+2}{\lambda_r} } \left( \iint_{{\cal Q}_{\rho}} u^r \ dx dt \right)^{2/\lambda_r} + \left(\frac{t}{\rho^2} \right)^{\frac{1}{1-m}}  \ .
\end{equation}
\end{proposition}

\begin{proof}
 Let $\sigma \in (0,1)$, $\rho>0$ and $t>0$ be fixed. For $ n \in \N_o$,  consider the increasing sequence of levels
\[k_n = k \left( 1 - \frac{1}{2^n} \right) \ , \]
being $k$ a positive number to be chosen along the proof satisfying
\begin{equation}\label{Ak}
k \ge \left(\frac{t}{\rho^2}\right)^{\frac{1}{(1-m)}}  \; 
\end{equation} 
the decreasing sequences of time levels and radii 
 \[t_n = t \left(\sigma + \frac{1-\sigma}{2^n} \right)  \ ,  \qquad
\rho_n= \rho \left(\sigma + \frac{1-\sigma}{2^n} \right)  \ , \]
from which we construct the sequences of nested and shrinking cubes and cylinders, respectively
 \[ \K_n= \prod_{i=1}^N \left\{ |x_i|< \left(\frac{t}{\rho^2}\right)^{\frac{m_i-m}{2(1-m)}} \rho_n \right\}  \quad \mathrm{and} \quad  \Q_n = \K_n \times [t-t_n,t] \ . \]
 
\vspace{.2cm}

\noindent Take smooth cutoff functions $\xi(x,t)= \xi_1(x) \xi_2(t)$ defined in $\Q_n$ and such that $\displaystyle{\xi_1(x)= \prod_{i=1}^N \xi_{1i}(x_i)}$ verifies
\[  \xi_{1i}=1  \ \mathrm{in} \ \K_{n+1}; \quad \xi_{1i}=0  \ \mathrm{in} \ \R^N \setminus \K_{n} \ ; \quad  |\xi_{1i}^{\prime}| \leq \left(\frac{t}{\rho^2}\right)^{\frac{m-m_i}{2(1-m)}} \frac{2^{n+1}}{(1-\sigma) \rho} \]  
and $\xi_2$, defined over  $[t-t_n, t]$, verifies
\[ \xi_2(\tau) =\left\{ 
\begin{array}{cl}
0 & , \ \tau \leq t-t_n \\[.8em]
1 & , \ t-t_{n+1} < \tau \leq t \ .
\end{array}\right.\]
In what follows we analyse each one of the two cases $1 \leq r \leq 2$ and $r>2$ separately. In both cases we consider that \eqref{condCrhoLrLinf} is not in force and take
\[ Y_n= \iint_{\Q_n} (u-k_{n})_+^r \, dxd\tau \, ;  \qquad   S= \sup_{\K_{\rho} \times [0, t]} u \qquad \mathrm{and} \qquad S_\sigma= \sup_{\K_{\sigma \rho} \times [\sigma t, t]} u \ . \]

\vspace{.2cm}

\noindent Case 1 : $1 \leq r \leq 2$

\vspace{.2cm}

\noindent Take test functions $\varphi= (u-k_{n+1})_+ \xi^2$ (which are admissible test functions due to the boundedness of $u$ and the smoothness of $\xi$), and consider the integration over $\Q_n$. Recalling the estimates for the cutoff functions, the structure conditions \eqref{SC1}-\eqref{SC2}-\eqref{SC3}, the fact that \eqref{condCrhoLrLinf} is not valid and applying Cauchy's inequality we arrive at
\[
    \sup_{t-t_n <\tau <t} \int_{\K_n} (u-k_{n+1})_+^2 \xi^2 \, dx\, + \frac{C_o}{2} m_1 \sum_{i=1}^N \iint_{\Q_n} u^{m_i-1} |u_{x_i}|^2 \xi^2 \ \chi[u>k_{n+1}] \, dxd\tau
 \]   
    \begin{eqnarray}\label{estkn+1}
    & \leq & \frac{2^{n+1}}{(1-\sigma)t}  \iint_{\Q_n} (u-k_{n+1})_+^2 \, dxd\tau \\ \nonumber 
    & & + \gamma (C_o,C_1) \frac{2^{2n}}{(1-\sigma)^2 \rho^2} \sum_{i=1}^N \left(\frac{t}{\rho^2}\right)^{\frac{m-m_i}{(1-m)}}\iint_{\Q_n} u^{m_i+1} \ \chi[u>k_{n+1}]\, \, dxd\tau \\\nonumber 
    & \leq &  \gamma (C_o,C_1) \frac{2^{4n}}{(1-\sigma)^2 t} \left\{1+ \sum_{i=1}^N \left(\frac{t}{\rho^2}\right)^{\frac{m-m_i}{(1-m)}+1} \frac{1}{k^{1-m_i}} \right\}\iint_{\Q_n} (u-k_{n})_+^2\, \, dxd\tau \\ \nonumber 
   & \leq & \gamma (C_o,C_1) \frac{2^{4n}}{(1-\sigma)^2 t} \left\{1+ \sum_{i=1}^N \left[\left(\frac{t}{\rho^2}\right)^{\frac{1}{(1-m)}} \frac{1}{k} \right]^{1-m_i} \right\}\iint_{\Q_n} (u-k_{n})_+^2 \, \, dxd\tau
\end{eqnarray}
by noting that
\begin{eqnarray*}
\iint_{\Q_n} (u-k_{n})_+^2 \, \, dxd\tau & \geq & \iint_{\Q_n} (u-k_{n})_+^2 \ \chi[u>k_{n+1}] \, \, dxd\tau = \iint_{\Q_n} u^2 \left(1- \frac{k_n}{u} \right)_+^2 \ \chi[u>k_{n+1}]\, \, dxd\tau\\
& = & \iint_{\Q_n} u^{1-m_i} \ u^{m_i+1} \left(1- \frac{k_n}{u} \right)_+^2 \ \chi[u>k_{n+1}]\, \, dxd\tau\\
& \geq & \left(\frac{k}{2} \right)^{1-m_i} \frac{1}{2^{2(n+1)}} \iint_{\Q_n} u^{m_i+1}\chi[u>k_{n+1}]\, \, dxd\tau \ .
\end{eqnarray*}
Now recall condition \eqref{Ak} and observe that
\begin{eqnarray*}
\sum_{i=1}^N \iint_{\Q_n} u^{m_i-1} |u_{x_i}|^2 \xi^2 \ \chi[u>k_{n+1}]  \, \, dxd\tau & \geq &\frac{1}{S} \ \sum_{i=1}^N \iint_{\Q_n} u^{m_i} |u_{x_i}|^2 \xi^2 \ \chi[u>k_{n+1}] \\
& = & \frac{4}{(m_N+2)^2 \ S}  \ \sum_{i=1}^N \iint_{\Q_n}  \left|\left((u-k_{n+1})_+^{(m_i+2)/2}\right)_{x_i}\right|^2 \xi^2\, \, dxd\tau \ . 
\end{eqnarray*}
From these estimates we get
\begin{equation*} \begin{aligned} \sup_{t-t_n \leq \tau \leq t} \int_{\K_n} (u-k_{n+1})_+^2 \xi^2\, dx &+ \frac{\gamma(C_o, m_1,m_N)}{S} \sum_{i=1}^N \iint_{\Q_n} \left|\left((u-k_{n+1})_+^{(m_i+2)/2}\right)_{x_i}\right|^2 \xi^2\, \, dxd\tau \\
& \leq \gamma (C_o,C_1,N) \frac{2^{4n}}{(1-\sigma)^2 t} \ S^{2-r} \ Y_n  \ .
 \end{aligned} \end{equation*}
Now consider 
\[\theta= \dfrac{N-2}{N} \ , \quad \alpha_i=\frac{m_i+2}{2} , \quad 2_{\alpha}^\star= 2^\star \dfrac{\sum_{i=1}^N \alpha_i}{N}, \quad \sigma= 2 \]
and \[q= \theta 2_{\alpha}^\star +2 (1-\theta) = m+2 + \frac{4}{N} >2 \geq r \ .\]
In what follows we derive an algebraic estimate involving the numbers $Y_n$ and $Y_{n+1}$: first we apply H\"{o}lder's inequality with exponent $q/r>1$ and then use the parabolic anisotropic embedding \eqref{PS} to obtain
\begin{eqnarray}\label{Yn}
 Y_{n+1} & \leq & \iint_{\Q_n}(u-k_{n+1})_+^r \xi^r\, \, dxd\tau \leq \left(\iint_{\Q_n}\left((u-k_{n+1})_+ \xi\right)^q \, dxd\tau \right)^{r/q} \ |\Q_n \cap [u>k_{n+1}]|^{1-r/q} \\ \nonumber
 & \leq & \gamma \left[\left(\sup_{t-t_n< \tau <t }\int_{\K_n}(u-k_{n+1})_+^2 \xi^2(x, \tau)\, dx\right)^{2/N} \ \prod_{i=1}^{N}\left(\iint_{\Q_n}| \left(\left((u-k_{n+1})_+ \xi \right)^{\alpha_{i}}\right)_{x_i}|^{2}\, dxd\tau\right)^{1/N} \right]^{r/q} \\[.8em] \nonumber
 & &  \hspace{4cm} \times |\Q_n \cap [u>k_{n+1}]|^{1-r/q} \\ [.8em] \nonumber
 & \leq & \gamma(N,C_o,C_1, m_1, m_N) \frac{b^n}{ ((1-\sigma)^2 t)^{\frac{(N+2)r}{Nq}}} \ k^{-\frac{r(q-r)}{q}} \ S^{\frac{((2-r)(N+2) +N)r}{Nq}} \ Y_n^{1+ \frac{2r}{Nq}} \ , \quad b= 2^{\frac{r(N(q-r)+4(N+2)}{Nq}} >1 \ , 
\end{eqnarray}
since
\[Y_n \geq \iint_{\Q_n} (u-k_{n})_+^r \chi[u>k_{n+1}] \ dx d\tau \geq \left(\frac{k}{2^{n+1}}\right)^r |\Q_n \cap [u>k_{n+1}]| \ . \]
By choosing $k>0$ 
\begin{equation}\label{choicek}
k=  \frac{\gamma}{((1-\sigma)^2 t)^{\frac{(N+2)}{N(q-r)}}} \ S^{\frac{(2-r)(N+2) +N}{N(q-r)}} \ \left( \iint_{\Q_o} u^r \ dx d\tau \right)^{\frac{2}{N(q-r)}} + \left(\frac{t}{\rho^2}\right)^{\frac{1}{(1-m)}}  \ , 
\end{equation}
not only one assures \eqref{Ak} but also one can apply Lemma \ref{fastgeomconv} to conclude that 
\begin{equation}\label{estSsigmaaniso}
S_\sigma \leq \frac{\gamma}{((1-\sigma)^2 t)^{\frac{(N+2)}{N(q-r)}}} \ S^{\frac{(2-r)(N+2) +N}{N(q-r)}} \ \left( \iint_{\Q_o} u^r \ dx dt \right)^{\frac{2}{N(q-r)}} + \left(\frac{t}{\rho^2}\right)^{\frac{1}{(1-m)}}  \ . 
\end{equation}
Define the sequences of real positive numbers 
\[\tilde{\rho}_o= \sigma \rho \ ; \qquad  \tilde{\rho}_n = \rho\left(\sigma + (1-\sigma) \sum_{j=1}^{n} 2^{-j}\right)\]
\[\tilde{t}_o= \sigma t \ ; \qquad  \tilde{t}_n = t\left(\sigma + (1-\sigma) \sum_{j=1}^{n} 2^{-j}\right)\]
and construct the cylinders
\[\tilde{\Q}_n = \tilde{\K}_n \times (t-\tilde{t}_n , t) \ , \qquad \tilde{\K}_n= \prod_{i=1}^N \left\{ |x_i|< \left(\frac{t}{\rho^2}\right)^{\frac{m_i-m}{2(1-m)}} \tilde{\rho}_n \right\}  \ . \]
Being $\displaystyle{S_n = \sup_{\tilde{\Q}_n} u}$, we then apply estimate \eqref{estSsigmaaniso} to the pair of cylinders $\tilde{\Q}_n$ and $\tilde{\Q}_{n+1}$ and use Young's inequality with exponent $\mu=\frac{N(q-r)}{(2-r)(N+2)+N}$ to get

\begin{eqnarray*}
S_n & \leq & \frac{\gamma}{((1-\sigma)^2 t)^{\frac{(N+2)}{N(q-r)}}}\ S_{n+1}^{\frac{(2-r)(N+2) +N}{N(q-r)}} \ \left( \iint_{\tilde{\Q}_{n+1}} u^r \ dx d\tau \right)^{\frac{2}{N(q-r)}} + \left(\frac{t}{\rho^2}\right)^{\frac{1}{(1-m)}}  \\
& \leq & \frac{1}{2} S_{n+1} + \frac{\gamma}{((1-\sigma)^2 t)^{\frac{N+2}{\lambda_r}}} \ \left( \iint_{\tilde{\Q}_{n+1}} u^r \ dx d\tau \right)^{\frac{2}{\lambda_r}} + \left(\frac{t}{\rho^2}\right)^{\frac{1}{(1-m)}}.
\end{eqnarray*}

\noindent By iteration (see Lemma \ref{iteration}) and taking $\sigma=1/2$ we finally get \eqref{estLrLinfty}.

\vspace{.5cm}

\noindent Case 2 : $r > 2$

\vspace{.2cm}

\noindent In this case we adopt a similar reasoning as the one presented before for $1 \leq r \leq 2$ but with a crucial difference: the boundedness of $u$ does not play a role as it did previously. We focus our attention in the main differences starting with the choice of the test functions: here we take $\varphi= (u-k_{n+1})_+^{r-1} \xi^2$. So, by considering test functions as such and integrating over $\Q_n$, recalling that \eqref{condCrhoLrLinf} is violated, one arrives at
\[\sup_{t-t_n \leq \tau \leq t} \int_{\K_n} (u-k_{n+1})_+^r \xi^2 (x,\tau) \ dx + \frac{ (r-1)r}{4} \ C_o m_1  \sum_{i=1}^N \iint_{\Q_n} u^{m_i-1}  (u-k_{n+1})_+^{r-2}|u_{x_i}|^2 \xi^2 \ dx d\tau\ 
 \]   
    \begin{eqnarray}\label{estknr}
    & \leq & \frac{2^{n+2}}{(1-\sigma)t}  \iint_{\Q_n} (u-k_{n+1})_+^r  \ dx d\tau\\ \nonumber 
    & & + \frac{\gamma (C_o,C_1,r)}{r-1} \frac{2^{2n}}{(1-\sigma)^2 \rho^2} \sum_{i=1}^N \left(\frac{t}{\rho^2}\right)^{\frac{m-m_i}{1-m}} \iint_{\Q_n} u^{m_i-1}  (u-k_{n+1})_+^{r} \ dx d\tau \\ \nonumber
    & & + \frac{\gamma (C_o,C_1,r)}{r-1} \frac{2^{2n}}{(1-\sigma)^2 \rho^2} \sum_{i=1}^N \left(\frac{t}{\rho^2}\right)^{\frac{m-m_i}{1-m}}
    \iint_{\Q_n} u^{m_i+1}  (u-k_{n+1})_+^{r-2} \ dx d\tau
    \\ \nonumber 
    & \leq &  \frac{\gamma (C_o,C_1,r)}{r-1} \frac{2^{2n}}{(1-\sigma)^2 t} \left\{ 1 + \sum_{i=1}^N \left[\left(\frac{t}{\rho^2}\right)^{\frac{1}{1-m}} \frac{1}{k} \right]^{1-m_i} \right\} \iint_{\Q_n} (u-k_{n+1})_+^{r} \ dx d\tau \\ \nonumber
    & & + \frac{\gamma (C_o,C_1,r)}{r-1} \frac{2^{2n}}{(1-\sigma)^2 t} \sum_{i=1}^N \left[\left(\frac{t}{\rho^2}\right)^{\frac{1}{1-m}} \frac{1}{k} \right]^{1-m_i} \iint_{\Q_n} u^2 (u-k_{n+1})_+^{r-2} \ dx d\tau \\ \nonumber
    & \leq & \frac{\gamma (C_o,C_1,r)}{r-1} \frac{2^{4n}}{(1-\sigma)^2 t}  \iint_{\Q_n} (u-k_{n+1})_+^{r}  \ dx d\tau \ .
\end{eqnarray}
These estimates have been obtained observing that
\[u > k_{n+1} > \frac{k}{2}\ \ \Rightarrow \ \ u^{m_i-1} \leq \frac{2}{k^{1-m_i}} \quad \forall i =1, \dots, N; \] 
\begin{eqnarray*} 
\iint_{\Q_n} (u-k_{n})_+^{r} \ dx d\tau  & \geq  &\iint_{\Q_n} (u-k_{n})_+^{r} \chi[u>k_{n+1}] \ dx d\tau \\
& = & 
\iint_{\Q_n} u^2 \left(1- \frac{k_n}{u}\right)_+^2  (u-k_{n})_+^{r-2} \ \chi[u>k_{n+1}]  \ dx dt \\
& \geq & \frac{1}{2^{2(n+1)}} \iint_{\Q_n} u^2 (u-k_{n+1})_+^{r-2} \ dx d\tau
\end{eqnarray*}
and taking $k>0$ as to satisfy \eqref{Ak}. As for the left-hand side of \eqref{estknr}, observe that
\begin{eqnarray*}
 \iint_{\Q_n} u^{m_i-1}  (u-k_{n+1})_+^{r-2}|u_{x_i}|^2 \xi^2 \ dx d\tau & \geq & \frac{1}{S} 
 \iint_{\Q_n}  (u-k_{n+1})_+^{m_i+r-2}|u_{x_i}|^2 \xi^2  \ dx d\tau   \\
& \geq &\frac{\gamma(r)}{S} \ \iint_{\Q_n}  \left|\left((u-k_{n+1})_+^{(m_i+r)/2}\right)_{x_i}\right|^2 \xi^2  \ dx d\tau   
\end{eqnarray*}  
and then 
\[\sup_{t-t_n \leq \tau \leq t} \int_{\K_n} (u-k_{n+1})_+^r \xi^2 (x,\tau) \ dx + \frac{\gamma(C_o, m_1, r)}{S} \ \sum_{i=1}^N \iint_{\Q_n}  \left|\left((u-k_{n+1})_+^{(m_i+r)/2}\right)_{x_i}\right|^2 \xi^2 \ dx d\tau \ 
 \] 
 \[ \leq   \gamma (C_o,C_1,r) \ \frac{2^{4n}}{(1-\sigma)^2 t} \ Y_n  \ .\]
Note that
\[ \hspace{-5cm}
 \iint_{\Q_n}  \left|\left(((u-k_{n+1})_+ \xi^2)^{(m_i+r)/2}\right)_{x_i}\right|^2 \ dx d\tau  \]
 \begin{eqnarray*}
 & \leq & 2 \iint_{\Q_n}  \left|\left((u-k_{n+1})_+^{(m_i+r)/2}\right)_{x_i}\right|^2 \xi^2\ dx dt + \frac{2^{2(n+2)}}{(1-\sigma)^2 \rho^2} (r+1)^2  \iint_{\Q_n} (u-k_{n+1})_+^{m_i+r}\, \ dx d\tau  \\[.8em]
& \leq & \frac{\gamma(r,C_o,C_1, m_1)}{(1-\sigma)^2 t} \ 2^{4n} \ S \left\{1 + \frac{t}{\rho^2} S^{m_i-1}\right\} Y_n \leq  \frac{\gamma(r,C_o,C_1, m_1)}{(1-\sigma)^2 t} \ 2^{4n} \ S \left\{1 + \frac{t}{\rho^2} \frac{1}{k^{1-m_i}} \right\} Y_n \ .
\end{eqnarray*}
The last inequality was obtained by recalling that  $ {k}/{2} < u \leq S$.
In what follows we derive an estimate for $Y_{n+1}$ by means of H\"{o}lder's inequality and the parabolic anisotropic embedding \eqref{PS}, now considering 
\[\theta= \dfrac{N-2}{N} \ , \quad \alpha_i=\frac{m_i+r}{2} , \quad 2_{\alpha}^\star= 2^\star \dfrac{\sum_{i=1}^N \alpha_i}{N} \ , \quad \sigma =r\ \ 
\]
and 
\[q= \theta 2_{\alpha}^\star +\sigma (1-\theta) = m+r +\frac{2r}{N} >r \ , \]
namely
\[Y_{n+1} \leq \frac{\gamma (C_o,C_1,m_1, m,r,N)}{\left((1-\sigma)^2 t\right)^{\frac{r(N+2)}{Nq}}} \ b^n \ S^{r/q} \ k^{-\frac{r(q-r)}{q}} \ Y_n^{1+\frac{2r}{Nq}} \ , \qquad b>1 \ .
\] 
The remainder of the proof is quite similar to the one presented for $1 \leq r \leq 2$, with the obvious changes: now, once we choose 
\begin{equation}\label{choicek2}
k=  \frac{\gamma}{((1-\sigma)^2 t)^{\frac{(N+2)}{N(q-r)}}} \ S^{\frac{1}{q-r}} \ \left( \iint_{\Q_o} u^r \ dx d\tau \right)^{\frac{2}{N(q-r)}} + \left(\frac{t}{\rho^2}\right)^{\frac{1}{(1-m)}}  \ , 
\end{equation}
we may conclude 
\begin{equation}\label{estSsigma}
S_\sigma \leq \frac{\gamma}{((1-\sigma)^2 t)^{\frac{(N+2)}{N(q-r)}}} \ S^{\frac{1}{q-r}} \ \left( \iint_{\Q_o} u^r \ dx d\tau \right)^{\frac{2}{N(q-r)}} + \left(\frac{t}{\rho^2}\right)^{\frac{1}{(1-m)}}  \ . 
\end{equation}
and afterward, by considering this estimate applied to the pair of cylinders $\tilde{\Q}_n$ and $\tilde{\Q}_{n+1}$, being 

\begin{eqnarray*}
S_n & \leq & \frac{\gamma}{((1-\sigma)^2 t)^{\frac{(N+2)}{N(q-r)}}}\ S_{n+1}^{\frac{(2-r)(N+2) +N}{N(q-r)}} \ \left( \iint_{\tilde{\Q}_{n+1}} u^r \ dx dt \right)^{\frac{2}{N(q-r)}} + \left(\frac{t}{\rho^2}\right)^{\frac{1}{(1-m)}}  \\
& \leq & \frac{1}{2} S_{n+1} + \frac{\gamma}{((1-\sigma)^2 t)^{\frac{N+2}{\lambda_r}}} \ \left( \iint_{\tilde{\Q}_{n+1}} u^r \ dx dt \right)^{\frac{2}{\lambda_r}} + \left(\frac{t}{\rho^2}\right)^{\frac{1}{(1-m)}}
\end{eqnarray*}
The proof is concluded once we choose $\sigma=1/2$ and use Lemma \ref{iteration}.
\end{proof}

\subsection*{Proving Theorem \ref{L1Linftest}}

\noindent Consider that \eqref{condCrho} is violated. From Proposition \ref{LrLinftest}, for $r=1$, and Theorem \ref{integralharnack} one obtains
\begin{eqnarray*}
    \sup_{\K{_{\rho/2}} \times [t/2,t]}u & \leq & \gamma \ t^{-\frac{N+2}{\lambda} } \left( \iint_{{\cal Q}_{\rho}} u \ dx d\tau \right)^{2/\lambda} + \left(\frac{t}{\rho^2} \right)^{\frac{1}{1-m}} \\
 & \leq & \gamma \ t^{-\frac{N}{\lambda} } \left(\sup_{0\leq \tau \leq t}  \int_{{\cal K}_{\rho}} u(x, \tau) \ dx  \right) ^{2/\lambda} + \left(\frac{t}{\rho^2} \right)^{\frac{1}{1-m}} \\
 & \leq & \gamma \ t^{-\frac{N}{\lambda} } \ \left\{ \inf_{0\leq\tau \leq t} \int_{{\cal K}_{2\rho}} u(x, \tau) \ dx + \left(\frac{t}{\rho^\lambda} \right)^{\frac{1}{1-m}} \right\}^{2/\lambda} + \left(\frac{t}{\rho^2} \right)^{\frac{1}{1-m}} \ , 
\end{eqnarray*}
which allows us to get \eqref{estL1Linfty}.

\section{$L^r_{loc}$ estimates backward in time} \label{S:Lrbackaniso}

\noindent This Section comprehends the presentation and proof of a local $L^r$ estimate backward in time: the following result states that the $L^r$-norm of (a locally bounded, nonnegative, local weak subsolution) $u$ to \eqref{gAPME} in an intrinsic anisotropic cube, located at any time level $0 \leq \tau \leq t$, can be bounded above by the $L^r$-norm of $u$ in a bigger cube located at the earliest time level (the bottom of a bigger (in space) cylinder).

\begin{proposition}{\bf [$L^r_{loc}$ estimates backward in time]}
\label{Lrbackest}\vskip0.1cm 
\noindent Let $u$ be a nonnegative, locally bounded, local weak subsolution to \eqref{gAPME}-\eqref{SC1}-\eqref{SC2}-\eqref{SC3} in $\Omega_T$. Assume $u \in L^r_{loc}(\Omega_T)$, for some $r>1$. There exists a positive constant $\gamma$, depending on $r, C_o, C_1, m_i$, such that, for all cylinders ${\Q}_{2\rho} \subset \Omega_T$, either
\begin{equation}\label{condCrho}
C \rho \left(\frac{t}{\rho^2}\right)^{\frac{m_i-m}{2(1-m)}} >1, \quad \mathrm{for} \ \mathrm{some} \ i=1, \cdots,N,
\end{equation}
or 
\begin{equation} \label{estLr}
\sup_{0 < \tau <t } \int_{\K_\rho} u^r (x, \tau) \ dx  \leq \gamma \left( \int_{{\K}_{2\rho}} u^r (x,0) \ dx   + \left(\frac{t}{\rho^{\lambda_r}} \right)^{\frac{1}{1-m}} \right)  \ .
\end{equation}
\end{proposition}

\begin{proof}
    Assume that \eqref{condCrho} fails. Let $\rho>0$, $t>0$ and $\sigma \in (0,1)$ be fixed and consider the cylinders $\Q_\rho \subset \Q_{(1+\sigma)\rho} \subset {\Q}_{2\rho} \subset \Omega_T$.  Take  $\varphi= f(u) \xi^2$, for
    \[f(u) = u^{r-1} \left( \frac{(u-k)_+}{u} \right)^s \ , \quad \max \{ r-1,1 \} < s< r \ , \quad k> 0\ , \]
    \[ F(u)= \int_k^u f(s) \ ds \ ,\]
and a time-independent smooth cutoff function $\xi \in C_o^{\infty}\left( \K_{(1+\sigma)\rho} \right)$ that equals one in $\K_\rho$, vanishes outside $\K_{(1+\sigma)\rho}$ and verifies $\displaystyle{|\xi_{x_i}| \leq \frac{1}{\sigma \rho} \left(\frac{t}{\rho^2}\right)^{\frac{m-m_i}{2(1-m)}}}$. 

\noindent Recalling \eqref{SC1}-\eqref{SC2}-\eqref{SC3}, the conditions on $\xi$, the failure of \eqref{condCrho}, the fact that
\[f^\prime (u) \geq (r-1) \frac{f(u)}{u} \qquad \mathrm{and} \qquad f(u) \leq u^{r-1} \ , \]
from the several terms appearing in the weak formulation \eqref{weaksol} one gets
\begin{eqnarray*}
  \sup_{0 < \tau <t } \int_{\K_\rho} F(u) (x,\tau) \ dx  & \leq & \int_{\K_{2\rho}} F(u) (x,0)  \ dx\\
  & & +  \frac{\gamma(C_o,C_1,r)}{(r-1) \sigma^2 \rho^2} \sum_{i=1}^N  \left\{ (C\rho)^2 +\left(\frac{t}{\rho^2}\right)^{\frac{m-m_i}{1-m}} + (C \rho)\left(\frac{t}{\rho^2}\right)^{\frac{m-m_i}{2(1-m)}} \right\} \\
  & & \hspace{2cm}
 \times \iint_{\Q_{(1+\sigma)\rho}} u^{m_i+r-1} \ \chi[u>k] \ dx d\tau \\
 & \leq & \int_{\K_{2\rho}} F(u) (x,0)  \ dx\\
& &
+ \frac{\gamma(C_o,C_1,r)}{(r-1) \sigma^2 \rho^2} \sum_{i=1}^N  \left(\frac{t}{\rho^2}\right)^{\frac{m-m_i}{1-m}}\iint_{\Q_{(1+\sigma)\rho}} u^{m_i+r-1} \ \chi[u>k] \ dx d\tau \ .
\end{eqnarray*}
In order to deal with the parabolic terms, on the one hand one proceeds as in \cite{FHV21} and chooses
\[ k^r \leq \frac{1}{2^{r+1} |\K_\rho|} \int_{\K_\rho} u^r (x, \tau) \ dx \]
to obtain the inferior bound to the left hand side
\[\sup_{0 < \tau <t } \int_{\K_\rho} F(u) (x,\tau) \ dx  \geq \gamma(r) \ \sup_{0 < \tau <t } \int_{\K_\rho} u^r (x,\tau) \ dx \ ; \]
on the other hand, since $F(u) \leq u^r$,
\[\int_{\K_{2\rho}} F(u) (x,0)  \ dx \leq \int_{\K_{2\rho}} u^r (x,0)  \ dx  \ .\]
As for the elliptic terms, we start by applying H\"older's inequality and then Young's inequality in each i-term to get
\[\frac{\gamma(C_o,C_1,r)}{(r-1) \sigma^2 \rho^2} \sum_{i=1}^N  \left(\frac{t}{\rho^2}\right)^{\frac{m-m_i}{1-m}}\iint_{\Q_{(1+\sigma)\rho}} u^{m_i+r-1} \ \chi[u>k] \ dx dt\]
\begin{eqnarray*}
    & \leq & \frac{\gamma(C_o,C_1,r)}{(r-1) \sigma^2 \rho^2} \sum_{i=1}^N  \left(\frac{t}{\rho^2}\right)^{\frac{m-m_i}{1-m}} \left(\iint_{\Q_{(1+\sigma)\rho}} u^r \ dx dt \right)^{\frac{m_i+r-1}{r}} \left( t \  2^{2N} \rho^N \right) ^{\frac{1-m_i}{r}} \\
    & \leq &  \sum_{i=1}^N  \left[ \frac{\gamma(C_o,C_1,r)}{(r-1) \sigma^2 }
    \left( \frac{t^r}{\rho^{\lambda_r}} \right)^{\frac{1-m_i}{(1-m)r}} \right] \
    \S_{\sigma}^{\frac{r-1+m_i}{r}} \ , \quad S_\sigma= \sup_{0 \leq \tau \leq t } \int_{\K_{(1+\sigma)\rho}} u^r (x,\tau) \ dx\\
   & \leq & \frac{1}{2} \S_\sigma + \gamma(C_o,C_1,r, \sigma,N,m_i) \ \left( \frac{t^r}{\rho^{\lambda_r}} \right)^{\frac{1}{1-m}} \ .
\end{eqnarray*}
Thereby
\begin{eqnarray*} 
\sup_{0 \leq \tau \leq t } \int_{\K_\rho} u^r (x,\tau) \ dx   & \leq & \gamma(r)  \int_{\K_{2\rho}} u^r (x,0)  \ dx \\
& & + \frac{1}{2} \S_\sigma + \gamma(C_o,C_1,r, \sigma,N,m_i) \ \left( \frac{t^r}{\rho^{\lambda_r}} \right)^{\frac{1}{1-m}} \ .
\end{eqnarray*}
Now consider
\[\rho_n= \rho \sum_{j=1}^n 2^{-j} \ ; \qquad \rho_{n+1} =(1+\sigma_n) \rho_n \ \ \Rightarrow \ \ \sigma_n \geq \frac{1}{2^{n+1}}\]
and
\[ S_n = \sup_{0 < \tau <t } \int_{\K_n} u^r (x,\tau) \ dx \ , \qquad \mathrm{being} \quad \K_n = \prod_{i=1}^N \left\{ |x_i|< \left(\frac{t}{\rho^2}\right)^{\frac{mi-m}{2(1-m)}}  \rho_n \right\} \  . \]
Therefore
\[ S_n \leq   \frac{1}{2} S_{n+1} +b^n \  \gamma(C_o,C_1,N, m_1, m_N, r) \left\{\int_{\K_{2\rho}} u^r (x,0)  \ dx   +  \left( \frac{t^r}{\rho^{\lambda_r}} \right)^{\frac{1}{1-m}} \right\} \ , \quad b>1 \ .\]
The result is now a consequence of Lemma \ref{iteration}.
\end{proof}

\begin{remark}
    The constant $\gamma$, determined along the proof, depends on $r$ in such a way that $\gamma\nearrow +\infty$ as $r\searrow 1$.
\end{remark}


\section{Extinction in finite time} \label{S:extinction}

In this Section we consider the Dirichlet problem \eqref{DP}-\eqref{SC1H}-\eqref{SC2H}, where $\Omega$ is a rectangular bounded domain in $\R^N$, $N>2$, and the initial data $0 \neq u_o \in L^\infty (\Omega)$. Let $u  \in C\left(\R^+; L^2(\Omega)\right)$, with  $u^{(m_i-1)/2} u_{x_i} \in L^2\left(\R^+;L^2(\Omega)\right), $ for all $i=1, \cdots, N$, be the unique, nonnegative, locally bounded solution to \eqref{DP}-\eqref{SC1H}-\eqref{SC2H}. In what follows we describe the procedure to determine the extinction of $u$ in a finite time $T^{\star}$, i.e. 
\[ u(x,t) =0 \ ,  \qquad \forall t\geq T^{\star}. \ \]
We start by multiplying the differential equation associated to \eqref{DP} by $\displaystyle{\varphi= \frac{u^r}{r+1}}$, for $r \geq m_i, \forall i$, and then we integrate over $\Omega$ to arrive at
\begin{eqnarray} \label{difequr+1}
0 & = & \frac{d}{dt} \int_\Omega u^{r+1} \ dx + C_o \frac{r}{r+1} \sum_{i=1}^N m_i \ \int_\Omega u^{m_i+r-2} |u_{x_i}|^2  \ dx\\ \nonumber
& \geq & \frac{d}{dt} \int_\Omega u^{r+1} \ dx + C_o m_1  \frac{4r}{(r+1)^3}   \sum_{i=1}^N  \ \int_\Omega \left| \left( u^{\alpha_i}\right)_{x_i} \right|^2 \ dx  \ ,  \qquad \alpha_i= \frac{m_i+r}{2}\\ \nonumber
& \geq & \frac{d}{dt} \int_\Omega u^{r+1} \ dx + \gamma C_o m_1  \frac{4r}{(r+1)^3} \  \left( \int_\Omega  u^{2_\alpha^\star} \ dx \right)^{\frac{2}{2^\star}} \ , \qquad \gamma=\gamma(N,m,r) \ ,
\end{eqnarray}
due to the elliptic anisotropic embedding \eqref{EE}. The next step consists in establishing a relation between the two integral terms $ \displaystyle{\int_\Omega u^{r+1} \ dx}$ and $\displaystyle{\int_\Omega  u^{2_\alpha^\star}   \ dx }$. For that purpose, we consider two possibilities: either the average $m$ is below or above the critical value $(N-2)/N$. We recall that when $m_i\equiv m$ for all $i=1, \dots, N$, this requirement distinguished between super and sub-critical exponents. Define 
    \[f(t) = \int_\Omega u^{r+1} (x,t) \ dx= ||u(\cdot, t)||_{L^{r+1}(\Omega)}^{r+1}  \] 
    and consider the alternative:

\vspace{.3cm}

\begin{enumerate} 

\item the average $m$ is in the sub-critical range, $\displaystyle{0<m< (N-2)/{N}}$:
\begin{enumerate}
    \item if $m_N < [N(1-m)]]/{2} -1$, take $r>m_N$ such that $r+1= [N(1-m)]/{2}$ and then $r+1=2_\alpha^\star$ and
    \[\int_\Omega u^{r+1} \ dx = \int_\Omega  u^{2_\alpha^\star} \ dx \]
    From \eqref{difequr+1} one obtains the ordinary differential inequality
     \begin{equation*}
        f^\prime + \gamma \ f^\beta \leq 0 \ , \quad \mathrm{where} \ \ \gamma=\gamma(C_o,m_1, N,m) , \ \ \beta= \frac{N-2}{N} <1
    \end{equation*}
    By integrating over the interval $[0,t]$, one gets
    \[ ||u(\cdot, t)||_{L^{r+1}(\Omega)} \leq ||u_o||_{L^{r+1}(\Omega)} \left( 1 - \frac{\gamma \ t}{||u_o||_{L^{r+1}(\Omega)}^{(1-\beta)(r+1)}}\right)_+^{\frac{1}{(1-\beta)(r+1)}}\]
    and
    \[0<T^\star \leq \gamma_1 \ ||u_o||_{L^{N(1-m)/2}(\Omega)}^{1-m} , \qquad \gamma_1=\gamma_1(C_o,m_1, N,m) \]
    
    \item if $m_N \geq  N(1-m)/2 -1$, take $r= m_N$. In this case  we have $r+1 < 2_\alpha^\star$ and 
    \[\int_\Omega u^{r+1} \ dx \leq \left( \int_\Omega  u^{2_\alpha^\star} \ dx \right)^{\frac{m_N+1}{2_\alpha^\star}} \ |\Omega|^{\frac{2_\alpha^\star -(1+m_N)}{2_\alpha^\star}} \ . \]
    Thereby
    \begin{equation*}
        f^\prime + \tilde{\gamma} \ |\Omega|^{(- 2_\alpha^\star +1+m_N)\frac{2}{2^\star(m_N+1)}} \ f^\beta \leq 0 \ , \quad \mathrm{where} \ \ \tilde{\gamma}=\tilde{\gamma}(C_o,m_1, m_N, N,m) ,  \, 
    \end{equation*} being $\beta= (m+m_N)/(1+m_N) <1$. Proceeding as before,
\[ ||u(\cdot, t)||_{L^{r+1}(\Omega)} \leq ||u_o||_{L^{r+1}(\Omega)} \left( 1 - \frac{\tilde{\gamma}  \ |\Omega|^{(- 2_\alpha^\star +1+m_N)\frac{2}{2^\star(m_N+1)}} \ t}{||u_o||_{L^{r+1}(\Omega)}^{(1-\beta)(r+1)}}\right)_+^{\frac{1}{(1-\beta)(r+1)}}\]
    and
    \[0<T^\star \leq \gamma_2  \ |\Omega|^{- \frac{N(m-1)+2 +2m_N}{N(m_N+1)^2}} \ ||u_o||_{L^{m_N+1}(\Omega)}^{1-m}  \ , \qquad \gamma_2 = \gamma_2(C_o,m_1, m_N, N,m) \ .\]
    
\end{enumerate}

\vspace{.3cm}

\item the average $m$ is in the supercritical range:  $\displaystyle{ (N-2)/{N} \leq m <1}$. 

In this case $ 2_\alpha^\star \geq r+1$, for all $r\geq m_i$. So we choose the smallest possible $r$, $r=m_N$,  and argue as in 1.(b) to get
\[ ||u(\cdot, t)||_{L^{r+1}(\Omega)} \leq ||u_o||_{L^{r+1}(\Omega)} \left( 1 - \frac{\tilde{\gamma}  \ |\Omega|^{(- 2_\alpha^\star +1+m_N)\frac{2}{2^\star(m_N+1)}} \ t}{||u_o||_{L^{r+1}(\Omega)}^{(1-\beta)(r+1)}}\right)_+^{\frac{1}{(1-\beta)(r+1)}}\]
    and
    \[0<T^\star \leq \gamma_2  \ |\Omega|^{- \frac{N(m-1)+2 +2m_N}{N(m_N+1)^2}} \ ||u_o||_{L^{m_N+1}(\Omega)}^{1-m}  \ , \qquad \gamma_2 = \gamma_2(C_o,m_1, m_N, N,m) \ .\]

    \end{enumerate}
\vspace{.3cm}

\noindent This proves Theorem \ref{extfinitetime}.

\vspace{.3cm}

\begin{remark}
    Observe that the assumption of boundedness of $u$ and $r \geq m_i$, for all $i=1, \cdots, N$, are necessary here consider admissible test functions of the form $\displaystyle{\varphi= {u^r}/({r+1})}$. 

\end{remark}

\vspace{.3cm}

\noindent As the proof shows, for $0<m< ({N-2})/{N}$ and $m_N < [{N(1-m)}]/{2} -1$, we did not make use of $\Omega$ being a bounded domain. Therefore, the extinction in finite time holds true also for the case of unbounded domains.

\begin{corollary}{\bf [Finite time of extinction in $\R^N$ ]}\label{extfinitetimewhole}
\vskip0.1cm 
\noindent Assume that  $u \in C(\R^+, L^2(\R^N))$, $u^{(m_i-1)/2} u_{x_i} \in L^2(\R^+, L^2(\R^N))$, $N>2$, is the unique nonnegative bounded solution to the Cauchy problem
\begin{equation} \label{CP}
\left\{
\begin{array} {ll}
u_t - \mathrm{div} A_o(x,t,u,Du) =0 & \ , \  (x,t) \in  \R^N\times \R^+,\\[.8em]
u(x,0)=u_o(x) \geq 0 & \ ,  \ x \in \R^N,
\end{array} \right.
\end{equation}
where $A_o$ satisfies \eqref{SC1H}-\eqref{SC2H} and $u_o \in L^\infty(\R^N) \cap L^1 (\R^N)$. Consider 
\[m_N <  \frac{N(1-m)}{2} -1 \quad \text{and} \quad 0<m<\frac{N-2}{N} \ . \]
Then there exits a finite positive time $T^\star$, depending on $C_o, N, m_1, m$, such that
\[u(x,t) \equiv 0 \ , \quad \mathrm{for} \ \ \mathrm{all} \ \ t \geq T^\star \ . \]
Moreover
\begin{equation}
    T^\star \leq  \gamma \ ||u_o||_{L^{N(1-m)/2}(\R^N)}^{1-m} 
\end{equation}
where the constant $\gamma$ only depends on $C_o, m_1, m, N$.
\end{corollary}

\begin{remark}
    Note that these two results, concerning the finite time of extinction, were obtained for a rectangular bounded space $\Omega$ and for the whole $\R^N$; no specification whatsoever in the choice of either an intrinsic or standard geometry.

    \noindent To the best of our knowledge, results concerning the existence and uniqueness of solutions to anisotropic Cauchy problems were obtained in a different setting than the one present here. In fact, in \cite{SongJian06} the authors considered
    \eqref{APME}, defined in $\R^N \times \R^+$, with initial data $u(x,0)=u_o(x) \in L^1(\R^N)$, and developed their work under the scope of 
    $L^1$-regular solutions: for each $\epsilon >0$, 
    $u \in C(\R^+_o; L^1(\R^N)) \cup C(\R^N \times \R^+) \cup L^{\infty}(\R^N \times [\epsilon, +\infty)) $.
    
\end{remark}

\vspace{.3cm}

\noindent The final part of this Section is devoted to the decay rate of extinction. For that purpose, let $\rho>0$ be fixed and take $t >{T^\star}/{2}$, being $T^\star$ the finite time of extinction given by Theorem \ref{extfinitetime} from which we construct the cylinder
\[ \K_{4\rho} \times (2t -T^\star, T^\star)=\prod_{i=1}^N \bigg{\{}|x_i|< \left( \frac{2(T^{\star}-t)}{\rho^2}\right) ^{(m_i-m)/(2(1-m))}4\rho \bigg{\}} \times (2t -T^\star, T^\star)\subset \Omega \times \R^+ \ .\]
From the integral Harnack-type inequality \eqref{estL1L1} one gets 
\[\sup_{2t -T^\star \leq \tau \leq  T^\star}  \int_{\K_{2\rho}} u(x, \tau) \ dx  \leq \gamma \  \left(\frac{T^\star -t}{\rho^\lambda} \right)^{\frac{1}{1-m}}  \ , \quad \lambda= N(m-1)+2 \]
and so, for all ${T^\star}/{2} < t <T^\star$, we have the following decay rate for the $L^1$-norm of $u$
\[  \int_{\K_{\rho}} u(x, t) \ dx  \leq \gamma \  \left(\frac{T^\star -t}{\rho^\lambda} \right)^{\frac{1}{1-m}} \ . \]
If in addition we consider $\lambda>0$, by recalling the local $L^1$-$L^\infty$ Harnack-type estimate \eqref{estL1Linfty}, we obtain Theorem \eqref{decayinfty}.




\section{The study of the anisotropic equation within the context of a standard geometry}\label{S:isotgeom}

\noindent The previously chosen geometry is closely related to the anisotropic behavior of our equation: considering radii $\rho_i=({t}/{\rho^2})^{(m_i-m)/(2(1-m))} \rho$ allowed us to obtain homogeneous estimates written for the arithmetic mean $m$. Therefore one can say that the anisotropic differential equation \eqref{gAPME} behaves, in this intrinsic anisotropic geometry, as the (general) porous medium type equation considered written for the arithmetic mean $m$. But $\rho_i$ either explodes or vanishes as soon as any asymptotic behavior (in time or space) needs to be checked. In  what follows we present similar results to the ones given in the previous Sections \ref{S:integralHarnack}, \ref{S:l1linftHarn}, \ref{S:Lrbackaniso} and \ref{S:extinction}, with the exception that now we consider the standard  geometry, meaning that we will consider the usual (isotropic) cubes
\[ K_\rho= \{ |x| < \rho \} \ , \qquad K_{a\rho}= \{ |x| < a\rho \}  \quad, \quad a>0,  \ \]
and the correspondent cylinders
\[ Q_\rho=K_\rho \times [0,t]  \ , \qquad  Q_{a\rho}= K_{a\rho}\times [0,t], \quad a>0 \ . \]
\noindent This will provide Harnack-type estimates that are stable from the geometric point of view: nevertheless, there is a price to be paid. In this Section, the differences between results obtained within the scope of the two adopted geometries will be clarified. 

\vspace{.2cm}

\noindent By presenting these two possible approaches, we provide the full picture; one may choose to consider the geometric setting that suits best the purpose to pursue. 

\subsection{Integral form of a Harnack-type inequality}

\begin{lemma}\label{lemma1isot}
\vskip0.1cm 
\noindent Let $u$ be a nonnegative, local weak supersolution to \eqref{gAPME}-\eqref{SC1}-\eqref{SC2}-\eqref{SC3} in $\Omega_T$. Let $\alpha \in \R$ be such that $0< \alpha +1 <1$, $0<\alpha +m_i <1$ and $0<m_i-\alpha <1$. There exists a positive constant $\gamma$, depending on $N, C_o, C_1,m_1, \alpha$, such that, for all cylinders $Q_{(1+\sigma)\rho} \subset \Omega_T$, for all $\sigma \in (0,1)$, either
\begin{equation}\label{condCrhoisotlemma}
C \rho >1 \ ,  
\end{equation}
or 
\begin{equation} \label{estL1isot}
\sum_{i=1}^N \iint_{Q_{\rho}} u^{m_i+\alpha-2} |u_{x_i}|^2  \ dx d\tau \leq \frac{\gamma}{\sigma^2} \ \left\{ 1+ \sum_{i=1}^N \left( \left(\frac{t}{\rho^2}\right)^{\frac{1}{1-m_i}} \frac{\rho^N}{S} \right)^{1-m_i} \right\} 
\S_\sigma^{1+\alpha} \ \rho^{-N\alpha}\ ,
\end{equation}
where $\displaystyle{S_\sigma = \sup_{0\leq \tau \leq t } \int_{ K_{(1+\sigma)\rho}} u(x,\tau) \ dx }$ and $ \ \displaystyle{S = \sup_{0 \leq \tau \leq t } \int_{{ K}_{\rho}} u(x,\tau) \ dx  }$.
\end{lemma}

\begin{proof}
    The proof follows closely the one of Lemma \ref{lemma1}, therefore we will just be focus on the changes to be considered. Assume \eqref{condCrhoisotlemma} does not hold and in the weak formulation \eqref{weaksol}, consider $\varphi=u^\alpha \xi(x)$, being $\xi(x)$
    a smooth cutoff function defined in ${K}_{(1+\sigma)\rho}$ and verifying $\xi=1$ in ${K}_{\rho}$, $\xi=0$ outside ${K}_{(1+\sigma)\rho}$ and $|D \xi| \leq {1}/{\sigma \rho} $, 
and integrate over $Q_{(1+\sigma)\rho}$. We then have 
\begin{eqnarray*}
0 & \leq & \int_{{ K}_{(1+\sigma)\rho}} u^{\alpha+1}(x,t) \ \xi(x) \ dx - \int_{{ K}_{(1+\sigma)\rho}} u^{\alpha+1}(x,0) \ \xi(x) \ dx - \iint_{{ Q}_{(1+\sigma)\rho}} u \left(u^{\alpha}\right)_t \ \xi \ dx dt \\
& & + \sum_{i=1}^N  \iint_{{Q}_{(1+\sigma)\rho}}A_i \left(\alpha \ u^{\alpha-1} \ u_{x_i} \ \xi + u^\alpha \ \xi_{x_i} \right) \ dx dt \\
& & -  \iint_{{ Q}_{(1+\sigma)\rho}} B \  u^\alpha \xi \ dx dt\\
& = & I_1 + I_2 + I_3.
\end{eqnarray*}
Observe that, as before, 
\begin{eqnarray*}
I_1 & \leq & \frac{2^{1+2N|\alpha|}}{\alpha+1} S_\sigma^{\alpha +1} \ \rho^{-N \alpha} \ ; 
\end{eqnarray*}
as for $I_2$ and $I_3$, remember that we are considering $C \rho \leq1$, $\alpha$ is negative and we have the estimate 
\[\iint_{{\cal Q}_{(1+\sigma)\rho}} u^{m_i+\alpha} \ dx d\tau \leq 
t \ S_\sigma ^{m_i + \alpha} \ (4\rho)^{N(1-(m_i+\alpha))} \ .
\]
Therefore 
\begin{eqnarray*}
I_2 & = & \alpha \sum_{i=1}^N  \iint_{Q_{(1+\sigma) \rho}} A_i u_{x_i}  \ u^{\alpha-1} \  \xi \ dx d\tau+ \sum_{i=1}^N  \iint_{Q_{(1+\sigma)\rho}} A_i \  u^\alpha \ \xi_{x_i} \ dx d\tau \\
& \leq & 2\alpha C_o \sum_{i=1}^N m_i  \iint_{Q_{(1+\sigma)\rho}} u^{m_i+\alpha-2} |u_{x_i}|^2 \xi^2 \ dx d\tau \\
& & + \frac{1}{|\alpha| (\sigma \rho)^2}  \sum_{i=1}^N  \left\{ (C\rho)^2 + \frac{C_1^2}{ C_o}   + C \rho  \right\} \iint_{Q_{(1+\sigma)\rho}} u^{m_i+\alpha} \ dx d\tau \\
& \leq & 2\alpha C_o \sum_{i=1}^N m_i  \iint_{{\cal Q}_{(1+\sigma)\rho}} u^{m_i+\alpha-2} |u_{x_i}|^2 \xi^2 \ dx d\tau \\
& & + \frac{\gamma(N, C_o,C_1)}{|\alpha| \sigma^2}  \sum_{i=1}^N  \left( 
\left(\frac{t}{\rho^2}\right)^{\frac{1}{1-m_i}}
\frac{\rho^N}{S} \right)^{1-m_i}  \ S_\sigma^{\alpha +1} \ \rho^{-N \alpha}   
\end{eqnarray*}
and 
\begin{eqnarray*}
I_3 & \leq &  - \alpha C_o  \sum_{i=1}^N  m_i \iint_{{\cal Q}_{1+\sigma}} u^{m_i-1+\alpha-1} |u_{x_i}|^2 \xi^2 \ dx d\tau \\
& & +  \frac{\gamma(N, C_o)}{|\alpha| \sigma^2} 
\sum_{i=1}^N  \left( 
\left(\frac{t}{\rho^2}\right)^{\frac{1}{1-m_i}}
\frac{\rho^N}{S} \right)^{1-m_i}  \ S_\sigma^{\alpha +1} \ \rho^{-N \alpha} \ .
\end{eqnarray*}
Finally, we get
\[\sum_{i=1}^N  \iint_{{\cal Q}_{(1+\sigma)\rho}} u^{m_i+\alpha-2} |u_{x_i}|^2 \xi^2 \ dx d\tau \leq \frac{\gamma(N,C_o, C_1, m_1)}{|\alpha|^2 (\alpha +1)} \frac{1}{\sigma^2} \left\{ 1+ \sum_{i=1}^N  \left( 
\left(\frac{t}{\rho^2}\right)^{\frac{1}{1-m_i}}
\frac{\rho^N}{S} \right)^{1-m_i}  \right\} \ S_\sigma^{\alpha +1} \ \rho^{-N \alpha} \  .
\]

\end{proof}

\vspace{.2cm}

\noindent This auxiliary lemma allows to derive an integral form of a Harnack-type inequality, namely

\begin{theorem}{\bf [Integral form of a Harnack-type inequality]}
\label{integralharnackisot}\vskip0.1cm 
\noindent Let $u$ be a nonnegative, local weak solution to \eqref{gAPME}-\eqref{SC1}-\eqref{SC2}-\eqref{SC3} in $\Omega_T$. There exists a positive constant $\gamma$, depending on $N, C_o, C_1, m_1$, such that, for all cylinders ${ Q}_{2\rho} \subset \Omega_T$, either
\begin{equation}\label{condCrhoisotth}
C \rho  >1 \ 
\end{equation}
or 
\begin{equation} \label{estL1isot}
\sup_{0<\tau < t}  \int_{{ K}_{\rho}} u(x, \tau) \ dx  \leq \gamma \ \left\{ \inf_{0<\tau < t} \int_{{ K}_{2\rho}} u(x, \tau) \ dx + \sum_{i=1}^N\left(\frac{t}{\rho^{\lambda_i}} \right)^{\frac{1}{1-m_i}} \right\} \ , 
\end{equation}
for $\lambda_i= N(m_i-1)+2$.  
\end{theorem}

\begin{proof}
    Assume that \eqref{condCrhoisotth} does not hold. Let $\rho>0$ and $t>0$ be fixed and construct the increasing sequence of (isotropic) cubes
\[K_{n}=  \left\{|x|< \rho_n \right\}\]
where 
\[\rho \leq \rho_n= \rho  \sum_{j=0}^n 2^{-j}  < 2 \rho \quad \mathrm{and} \quad \rho_{n+1}=(1+\sigma_n) \rho_n \ \Rightarrow \ \sigma_n \geq \frac{1}{2^{n+1}} \ .\]
In the weak formulation take $\varphi= \xi(x)$, being $\xi(x)$ a smooth cutoff function that: equals $1$ in $K_{n}$, vanishes outside $K_{n+1}$, verifies $|D \xi | \leq {2^{n+1}}/{\rho} $, and consider the integration over $K_{n+1} \times [\tau_1, \tau_2] \subset K_{2 \rho} \times [0,t]$. We then get, 
\begin{eqnarray} \label{firstestintharnackisot}
\int_{K_{n}}  u(x, \tau_1) \, dx & \leq  &\int_{K_{n+1}} u(x, \tau_2)  \, dx \nonumber \\ 
& & +C_1 \frac{2^{n+1}}{\rho}  \sum_{i=1}^N m_i \int_{\tau_1}^{\tau_2}\int_{K_{n+1}}  u^{m_i-1} |u_{x_i} | \, dxd\tau  + C \sum_{i=1}^N m_i \int_{\tau_1}^{\tau_2}\int_{K_{n+1}}u^{m_i-1} |u_{x_i} |\, dxd\tau\nonumber \\
& & + C \frac{2^{n+1}}{\rho}  \sum_{i=1}^N  
\int_{\tau_1}^{\tau_2}\int_{K_{n+1}}  u^{m_i}  \, dxd\tau  + C^2 \sum_{i=1}^N \int_{\tau_1}^{\tau_2}\int_{K_{n+1}} u^{m_i} \, dxd\tau  \ .
\end{eqnarray} 
Now choose $\tau_2 \in [0,t]$ such that
\[ \int_{K_{2\rho}} u (x, \tau_2) \ d\tau = 
\inf_{0\leq \tau \leq t}  \int_{K_{2\rho}} u(x,\tau) \ dx = \I \]
 and denote  $\displaystyle{S_n=\sup_{0\leq \tau\leq t}  \int_{K_{n}} u(x, \tau) \ dx }$ and $Q_{n}= K_n \times [0,t]$. Taking this into account, recalling $C \rho \leq 1$ and Lemma \ref{lemma1isot}, from the previous inequality \eqref{firstestintharnackisot} we get
 \begin{eqnarray*}
     S_n & \leq & \I +  \frac{2^{n+1}}{\rho} (C_1 + C\rho) \sum_{i=1}^N   \iint_{Q_{n+1}}  u^{m_i-1} |u_{x_i} | \, dxd\tau \\
& & + \frac{2^{n+1}}{\rho^2} (C \rho + (C\rho)^2 ) \sum_{i=1}^N  \iint_{Q_{n+1}}  u^{m_i}  \, dxd\tau  \\
& \leq & \I + \frac{2^{2n}}{\rho} \gamma (C_o,C_1,N,m_1)  \sum_{i=1}^N \left\{ 1 + \sum_{j=1}^N \left( \frac{t}{\rho^2} \right) \left( \frac{\rho^N}{S_o} \right)^{1-m_j} \right\}^{1/2} \ 
\sqrt{t} \ S_{n+1}^{\frac{1+m_i}{2}} \ \rho^{{\frac{N(1-m_i)}{2}}} \\
&  &  + \frac{2^n}{\rho^2} \gamma(N)\  t \sum_{i=1}^N S_{n+1}^{1+m_i} \ \rho^{N(1-m_i)} \\
& \leq & \I + \sum_{i=1}^N   \ S_{n+1}^{\frac{1+m_i}{2}} \left[ 2^{2n}    \gamma (C_o,C_1,N,m_1)   \ \rho^{{\frac{N(1-m_i)}{2}}} \left( \frac{t}{{\rho^2}}\right)^{1/2} \right] \\
 &  &  + \sum_{i=1}^N S_{n+1}^{m_i} \ \left[ \ \gamma(N) 2^{n} \ \rho^{N(1-m_i)} \frac{t}{\rho^2} \right] \ .
 \end{eqnarray*}
The last estimate was obtained assuming, without loss of generality (otherwise there is nothing more to be done), the inequality
\begin{equation} \label{kondizionen} S= S_o > \sum_{i=1}^N \left( \frac{t}{\rho^{N(m_i-1)+2}} \right)^{\frac{1}{1-m_i}}  \ . \end{equation}
Our next step will be to apply in each i-term of the second and third terms, Young's inequality with $\epsilon_i= {\epsilon}/{2N}$, for small $\epsilon> 0$, with exponents $\mu_i={2}/({1-m_i})$ and $\nu_i= {1}/({1-m_i})$, respectively, and then get
\[ S_n \leq  \epsilon \ S_{n+1} + b^n \ \gamma(N,C_o,C_1,m_N,\epsilon) \left[ \sum_{i=1}^N \left( \frac{t}{\rho^{N(m_i-1)+2}} \right)^{\frac{1}{1-m_i}}   + \I \right]\ , \quad b=2^{\frac{1}{1-m_N}}>1 \ . \]
The proof is concluded once we choose $\epsilon \in (0,1)$ such that $\epsilon b =1/2$ and let $n \rightarrow \infty$.
\end{proof}

\subsection{ Local $L^1$-$L^\infty$ Harnack-type estimates and Decay rate of extinction}

\begin{proposition}{\bf [$L^r_{loc}-L^{\infty}_{loc}$ estimates]}
\label{LrLinftestisot}\vskip0.1cm 
\noindent Let $u$ be a nonnegative, locally bounded,  local weak subsolution to \eqref{gAPME}-\eqref{SC1}-\eqref{SC2}-\eqref{SC3} in $\Omega_T$. Let $r \geq 1 $ be such that $\lambda_r= N(m-1)+2r >0$. There exists a positive constant $\gamma$, depending on $N, C_o, C_1, m_i$, such that, for all cylinders $Q_\rho = K_\rho \times [0, t] \subset \Omega_T$, either
\begin{equation}\label{condCrhoisot}
C \rho >1 \ , 
\end{equation}
or 
\begin{equation} \label{estLrisot}
\sup_{K{_{\rho/2}} \times [t/2,t]}u \leq \gamma \ t^{-\frac{N+2}{\lambda_r} } \left( \iint_{Q_\rho} u^r \ dx d\tau \right)^{2/\lambda_r} + \sum_{i=1}^N\left(\frac{t}{\rho^2} \right)^{\frac{1}{1-m_i}}  \ 
\end{equation}
\end{proposition}

\begin{proof}
Let $\sigma \in (0,1)$, $\rho>0$ and $t>0$ be fixed and for $ n=0, 1, \cdots$,  consider the decreasing sequences of time levels
 \[t_n = t \left(\sigma + \frac{1-\sigma}{2^n} \right)\]
 and of radii
 \[\rho_n= \rho \left(\sigma + \frac{1-\sigma}{2^n} \right)  \ , \]
 from which we construct the sequences of nested and shrinking cubes and cylinders, respectively,
 \[ K_n= \left\{ |x|< \rho_n \right\}  \quad \mathrm{and} \quad  Q_n = K_n \times [t-t_n,t] \ . \]
Consider the increasing sequence of levels, being $k$ a positive to be fixed, 
\[k_n = k \left( 1 - \frac{1}{2^n} \right) \quad n=0,1,\cdots \ . \]

\vspace{.2cm}

\noindent Take smooth cutoff functions $\xi(x,t)= \xi_1(x) \xi_2(t)$ defined in $Q_n$ and such that $\xi_1$ verifies
\[  \xi_{1}=1  \ \mathrm{in} \ K_{n+1}; \quad \xi_{1}=0  \ \mathrm{in} \ \R^N \setminus K_{n} \ ; \quad  |D\xi_{1}| \leq \frac{2^{n+1}}{(1-\sigma) \rho} \]  
and $\xi_2$, defined over the interval $[t-t_n, t]$, verifies
\[ \xi_2(\tau) =\left\{ 
\begin{array}{cl}
0 & , \ \tau \leq t-t_n \\[.8em]
1 & , \ t-t_{n+1} < \tau \leq t \ .
\end{array}\right.\]
assume \eqref{condCrhoisot} is not in force and define
\[ Y_n= \iint_{Q_n} (u-k_{n}\ )_+^r \ dx d\tau\ ;  \qquad   S= \sup_{K_{\rho} \times [0, t]} u \qquad \mathrm{and} \qquad S_\sigma= \sup_{K_{\sigma \rho} \times [\sigma t, t]} u \ . \]

\vspace{.2cm}

\noindent Consider $1 \leq r \leq 2$. Arguing as in the correspondent part of the proof of Proposition \ref{LrLinftest}, we arrive at
\[
    \sup_{t-t_n \leq \tau \leq t} \int_{K_n} (u-k_{n+1})_+^2 \xi^2 \ dx + \frac{C_o}{2} m_1 \sum_{i=1}^N \iint_{Q_n} u^{m_i-1} |u_{x_i}|^2 \xi^2 \ \chi[u>k_{n+1}]  \ dx d\tau
 \]   
    \begin{eqnarray}\label{estkn+1}
    & \leq & \frac{2^{n+2}}{(1-\sigma)t}  \iint_{Q_n} (u-k_{n+1})_+^2  \ dx d\tau\\ \nonumber 
    & & + \gamma (C_o,C_1) \frac{2^{2n}}{(1-\sigma)^2 \rho^2} \sum_{i=1}^N \iint_{Q_n} u^{m_i+1} \ \chi[u>k_{n+1}] \ dx d\tau \\\nonumber 
    & \leq &  \gamma (C_o,C_1) \frac{2^{4n}}{(1-\sigma)^2 t} \left\{1+ \sum_{i=1}^N \left(\frac{t}{\rho^2}\right) \frac{1}{k^{1-m_i}} \right\}\iint_{\Q_n} (u-k_{n})_+^2 \ dx d\tau \ .
\end{eqnarray}
For the sake of obtaining an estimate for the right-hand side that is independent of the index $i$, we consider a level $k$ such that
\begin{equation}\label{Aki}
k \ge \bigg( \frac{t}{rho^2} \bigg)^{\frac{1}{1-m_i}}\ , \qquad \mathrm{for} \ \ \mathrm{all} \ \ i=1, \cdots, N \ .
\end{equation}
This is the price to be paid for having a homogeneous right-hand side, using the nonhomogeneous assumption \eqref{Aki}. As for the left-hand side we get the inferior bound 
\[ \sum_{i=1}^N  \iint_{Q_n} u^{m_i-1} |u_{x_i}|^2 \xi^2 \ \chi[u>k_{n+1}]  \ dx d\tau  \geq \frac{4}{9 S} \sum_{i=1}^N \iint_{Q_n}  \left|\left((u-k_{n+1})_+^{(m_i+2)/2}\right)_{x_i}\right|^2 \xi^2 \ dx d\tau \ , \]
and by combining both estimates 
\[\sup_{t-t_n \leq \tau \leq t} \int_{K_n} (u-k_{n+1})_+^2 \xi^2  \ dx + \frac{\gamma(C_o,m_1)}{S} \sum_{i=1}^N \iint_{Q_n} \left|\left((u-k_{n+1})_+^{(m_i+2)/2}\right)_{x_i}\right|^2 \xi^2  \ dx d\tau \]
\[\leq \gamma (C_o,C_1,N) \frac{2^{4n}}{(1-\sigma)^2 t} \ S^{2-r} \ Y_n  \ . \]   
We now consider 
\[q= \theta 2_{\alpha}^\star +2 (1-\theta), \quad \theta= \dfrac{N-2}{N} \ , \quad \alpha_i=\frac{m_i+2}{2} , \quad 2_{\alpha}^\star= 2^\star \dfrac{\sum_{i=1}^N \alpha_i}{N} \ , \]
apply H\"{o}lder's inequality with exponent $q/r>1$ and then use the anisotropic embedding \eqref{PS} to obtain
\[ Y_{n+1} \leq  \gamma(N,C_o, C_1, m_1) \frac{b^n}{ ((1-\sigma)^2 t)^{\frac{(N+2)r}{Nq}}} \ k^{-\frac{r(q-r)}{q}} \ S^{\frac{((2-r)(N+2) +N)r}{Nq}} \ Y_n^{1+ \frac{2r}{Nq}} \ , \quad b>1 \ . \]

\noindent To use Lemma \ref{fastgeomconv} and accommodate assumption \eqref{Aki}, we take 
\begin{equation}\label{choiceki}
k=  \frac{\gamma}{((1-\sigma)^2 t)^{\frac{(N+2)}{N(q-r)}}} \ S^{\frac{(2-r)(N+2) +N}{N(q-r)}} \ \left( \iint_{Q_o} u^r \ dx d\tau \right)^{\frac{2}{N(q-r)}} + \sum_{i=1}^N\left(\frac{t}{\rho^2}\right)^{\frac{1}{1-m_i}}  \ , 
\end{equation}
and then get 
\begin{equation}\label{estSsigmai}
S_\sigma \leq \frac{\gamma}{((1-\sigma)^2 t)^{\frac{(N+2)}{N(q-r)}}} \ S^{\frac{(2-r)(N+2) +N}{N(q-r)}} \ \left( \iint_{Q_o} u^r \ dx d\tau \right)^{\frac{2}{N(q-r)}} + \sum_{i=1}^N\left(\frac{t}{\rho^2}\right)^{\frac{1}{1-m_i}}  \ . 
\end{equation}
From this point on there are so significant changes to be made to the proof presented within the anisotropic geometry setting.

\noindent Let us now consider $r>2$. This case follows quite closely the reasoning presented before; we only present the main differences $\varphi= (u-k_{n+1})_+^{r-1} \xi^2$. So by considering test functions as such and integrating over $Q_n$ one arrives at
\[\sup_{t-t_n \leq \tau \leq t} \int_{K_n} (u-k_{n+1})_+^r \xi^2 (x,\tau) \ dx + \frac{ (r-1)r}{4} \ C_o m_1  \sum_{i=1}^N \iint_{Q_n} u^{m_i-1}  (u-k_{n+1})_+^{r-2}|u_{x_i}|^2 \xi^2  \ dx d\tau \ 
 \]   
    \begin{eqnarray}\label{estknrisot}
    & \leq & \frac{2^{n+2}}{(1-\sigma)t}  \iint_{Q_n} (u-k_{n+1})_+^r  \ dx d\tau \nonumber \\
    & & + \frac{\gamma (C_o,C_1,r)}{r-1} \frac{2^{2n}}{(1-\sigma)^2 \rho^2} \sum_{i=1}^N \iint_{Q_n} u^{m_i-1}  (u-k_{n+1})_+^{r} \ dx d\tau \nonumber  \\
    & & + r(r-1) C^2 \sum_{i=1}^N   \iint_{Q_n} u^{m_i+1}  (u-k_{n+1})_+^{r-2} \ dx d\tau \nonumber \\
    & & + C^2 \sum_{i=1}^N   \iint_{Q_n} u^{m_i}  (u-k_{n+1})_+^{r-1} \ dx d\tau \nonumber \\ 
    & \leq & \gamma (C_o,C_1,r)\frac{2^{2n}}{(1-\sigma)^2 t}  \iint_{Q_n} (u-k_{n+1})_+^{r}  \ dx d\tau
\end{eqnarray}
These estimates were obtained, considering $ C\rho \leq 1$, assuming $k$ verifies \eqref{Aki}, noting that
\[u >k_{n+1} > \frac{k}{2} \ \Rightarrow \ u^{m_i-1} \leq \frac{2}{k^{1-m_i}} \quad \forall i=1, \dots, N ; \] and using the same reasoning as before (in the anisotropic setting) to estimate the several integrals in terms of the $L^r$-norm of the truncated functions $(u-k_n)_+$. As for the left hand side, once again one has
\[ \iint_{Q_n} u^{m_i-1}  (u-k_{n+1})_+^{r-2}|u_{x_i}|^2 \xi^2 \ dx d\tau \geq \frac{\gamma(r)}{S} \ \iint_{Q_n}  \left|\left((u-k_{n+1})_+^{(m_i+r)/2}\right)_{x_i}\right|^2 \xi^2 \ dx d\tau\  \] 
thereby 
\[\sup_{t-t_n \leq \tau \leq t} \int_{K_n} (u-k_{n+1})_+^r \xi^2 (x,\tau) \ dx + \frac{\gamma(C_o, m_1, r)}{S} \ \sum_{i=1}^N \iint_{Q_n}  \left|\left((u-k_{n+1})_+^{(m_i+r)/2}\right)_{x_i}\right|^2 \xi^2 \ dx d\tau \] 
 \[  \leq   \gamma (C_o,C_1,r) \frac{2^{2n}}{(1-\sigma)^2 t} \ Y_n  \ . \]
Consider \[q= \theta 2_{\alpha}^\star + r (1-\theta) \ , \quad \theta= \dfrac{N-2}{N} \ , \quad \alpha_i=\frac{m_i+r}{2} , \quad 2_{\alpha}^\star= 2^\star \dfrac{\sum_{i=1}^N \alpha_i}{N} \ , \] and use H\"{o}lder's inequality and the parabolic anisotropic embedding \eqref{PS}, to get 
\[Y_{n+1} \leq \frac{\gamma (C_o,C_1,m_1, m,r,N)}{\left((1-\sigma)^2 t\right)^{\frac{r(N+2)}{Nq}}} \ b^n \ S^{r/q} \ k^{-\frac{r(q-r)}{q}} \ Y_n^{1+\frac{2r}{Nq}}
\]
The remaining of the proof is quite similar to the one presented for $1 \leq r \leq 2$, with the obvious changes: now, once we choose 
\begin{equation}\label{choicek2}
k=  \frac{\gamma}{((1-\sigma)^2 t)^{\frac{(N+2)}{N(q-r)}}} \ S^{\frac{1}{q-r}} \ \left( \iint_{Q_o} u^r \ dx d\tau \right)^{\frac{2}{N(q-r)}} + \sum_{i=1}^N\left(\frac{t}{\rho^2}\right)^{\frac{1}{1-m_i}}  \ , 
\end{equation}
we may conclude 
\begin{equation}\label{estSsigma}
S_\sigma \leq \frac{\gamma}{((1-\sigma)^2 t)^{\frac{(N+2)}{N(q-r)}}} \ S^{\frac{1}{q-r}} \ \left( \iint_{Q_o} u^r \ dx d\tau \right)^{\frac{2}{N(q-r)}} + \sum_{i=1}^N\left(\frac{t}{\rho^2}\right)^{\frac{1}{1-m_i}}  \ . 
\end{equation}
and afterward, by considering this estimate applied to the pair of cylinders $\tilde{Q}_n$ and $\tilde{Q}_{n+1}$ (for $\tilde{\rho}_n$ and $\tilde{t}_n$ defined as in the anisotropic geometry case studied in Section \ref{S:l1linftHarn}) 
\begin{eqnarray*}
S_n & \leq & \frac{\gamma}{((1-\sigma)^2 t)^{\frac{(N+2)}{N(q-r)}}}\ S_{n+1}^{\frac{1}{q-r}} \ \left( \iint_{\tilde{Q}_{n+1}} u^r \ dx d\tau \right)^{\frac{2}{N(q-r)}} + \sum_{i=1}^N \left(\frac{t}{\rho^2}\right)^{\frac{1}{1-m_i}}  \\
& \leq & \frac{1}{2} S_{n+1} + \frac{\gamma}{((1-\sigma)^2 t)^{\frac{N+2}{\lambda_r}}} \ \left( \iint_{\tilde{Q}_{n+1}} u^r \ dx d\tau \right)^{\frac{2}{\lambda_r}} + \sum_{i=1}^N \left(\frac{t}{\rho^2}\right)^{\frac{1}{1-m_i}}
\end{eqnarray*}
and the result follows by applying Lemma \ref{iteration}.

\end{proof}

\noindent By combining Proposition \ref{LrLinftestisot}, for $r=1$, and Theorem \ref{integralharnackisot} one obtains

\begin{theorem} {\bf [$L^1_{loc}-L^{\infty}_{loc}$ Harnack-type estimate]}
\label{L1Linftestisot}\vskip0.1cm 
\noindent Let $u$ be a nonnegative, locally bounded, local weak solution to \eqref{gAPME}-\eqref{SC1}-\eqref{SC2}-\eqref{SC3} in $\Omega_T$, for $m$ in the supercritical range, {\it{i.e.}} $\displaystyle{m >(N-2)/N}$. There exists a positive constant $\gamma$, depending on $N, C_o, C_1, m_i$, such that, for all cylinders $Q_{2\rho} \subset \Omega_T$, either
\begin{equation}\label{condCrhoisotl1linf}
C \rho  >1 \ , 
\end{equation}
or 
\begin{equation} \label{estL1Linfisot}
\sup_{K{_{\rho/2}} \times [t/2,t]}u  \leq  \gamma \ t^{-\frac{N}{\lambda} } \left(\inf_{ 0 \leq  \tau \leq \rho} \int_{K_{2\rho}} u (x, \tau) \ dx \right)^{2/\lambda} + \gamma \sum_{i=1}^N \left(\frac{t}{\rho^2} \right)^{\frac{1}{1-m_i}} +  \gamma \sum_{i=1}^N  \left(\frac{t}{\rho^2} \right)^{\frac{\lambda_i}{\lambda(1-m_i)}} \ ,
\end{equation}
where $\lambda_i= N(m_i-1)+2$ and $\lambda= N(m-1)+2$.
\end{theorem}

\vspace{.3cm}

\begin{remark}
 Note that the $L^1_{loc}$-$L^{\infty}_{loc}$ Harnack-type estimate was derived just asking for $m>(N-2)/{N}$, that is, considering $\lambda>0$. However, with no further conditions on the exponents $m_i$, $i=1,\dots, N$, the exponents $\lambda_i$ do not have a constant sign. The exponents' positivity  is a relevant feature when deriving  a decay rate of extinction.    
\end{remark}

\vspace{.2cm}

\noindent In what follows we present the decay rate of extinction in the case of the standard geometry. We start by considering $\rho>0$ fixed, $T^\star/2<t<T^\star$, being $T^\star$ the finite time of extinction given by \eqref{Tstar}, and then construct the cylinder
\[ K_{4\rho} \times [2t -T^\star, T^\star]\subset \Omega \times \R^+ \ .\]
From the integral Harnack-type inequality \eqref{estL1isot} 
\[\sup_{2t -T^\star<\tau < T^\star}  \int_{K_{\rho}} u(x, \tau) \ dx  \leq \gamma \  \sum_{i=1}^N\left(\frac{T^\star -t}{\rho^{\lambda_i}} \right)^{\frac{1}{1-m_i}}    \]
and so, for all $\dfrac{T^\star}{2} < t <T^\star$, we have the following decay rate for the $L^1$-norm of $u$
\[  \int_{K_{\rho}} u(x, t) \ dx  \leq \gamma \  \sum_{i=1}^N\left(\frac{T^\star -t}{\rho^{\lambda_i}} \right)^{\frac{1}{1-m_i}} \ . \]

\noindent Observe that if furthermore we assume $\lambda >0$, by recalling the local $L^1$-$L^\infty$ Harnack-type estimate \eqref{estL1Linfisot}, we have
\begin{equation}\label{decaystandardgeom}
||u(\cdot, t)||_{L^{\infty}(K_\rho)}  \leq  \gamma \ \sum_{i=1}^N \left(\frac{T^\star- t}{\rho^2} \right)^{\frac{1}{1-m_i}} +  \gamma \sum_{i=1}^N  \left(\frac{T^\star-t}{\rho^2} \right)^{\frac{\lambda_i}{\lambda(1-m_i)}} \ , \quad \dfrac{T^\star}{2} < t <T^\star \ .
\end{equation}

\noindent By a simple analysis of each one of the two terms of \eqref{decaystandardgeom}, and being forced to assume a stronger condition on the exponents $m_i$ (for instance observe that $\lambda_i$ may be of either sign at this point), we derive the following decay rate in the standard geometry.

\begin{theorem}\label{decayinftyisot}
{\bf [Decay rate of extinction]}
\vskip0.1cm 
\noindent
In the setting of Theorem \ref{extfinitetime}, consider the smallest exponent $m_1>(N-2)/N$ to be supercritical. Then, there exists a positive constant $\gamma$, depending upon $C_o,m_1, m_N, N,m$, such that, for all ${T^\star}/{2} < t <T^\star$, 
\begin{enumerate}
\item[(a)] when $(T^\star- t)/\rho^2 \leq 1$, 
\[ ||u(\cdot, t)||_{L^{\infty}(K_\rho)}  \leq \gamma \  \left(\frac{T^\star -t}{\rho^2} \right)^{\frac{\lambda_1}{\lambda(1-m_1)}} \ \ ;   \]
\item[(b)] if otherwise $(T^\star- t)/\rho^2 \geq 1$, 
\[ ||u(\cdot, t)||_{L^{\infty}(K_\rho)}  \leq \gamma \  \left(\frac{T^\star -t}{\rho^2} \right)^{\frac{\lambda_N}{\lambda(1-m_N)}} \ \ .   \]
\end{enumerate}
\end{theorem}

\vspace{.3cm}
\noindent We conclude with the study of backward estimates for the case of the standard geometry.
\subsection{Local $L^r$ estimates backward in time }

\begin{proposition}{\bf [$L^r_{loc}$ estimates backward in time]}
\label{Lrbackestisot}\vskip0.1cm 
\noindent Let $u$ be a nonnegative, locally bounded, local weak subsolution to \eqref{gAPME}-\eqref{SC1}-\eqref{SC2}-\eqref{SC3} in $\Omega_T$. Assume $u \in L^r_{loc}(\Omega_T)$, for some $r>1$. There exists a positive constant $\gamma$, depending on $r, C_o, C_1, m_i$, such that, for all cylinders ${Q}_{2\rho} \subset \Omega_T$, either
\begin{equation}\label{condCrhoisotLr}
C \rho  >1 \ ,
\end{equation}
or 
\begin{equation} \label{estLrisot}
\sup_{0 < \tau <t } \int_{K_\rho} u^r (x, \tau) \ dx  \leq \gamma \left( \int_{{K}_{2\rho}} u^r (x,0) \ dx   + \sum_{i=1}^N \left(\frac{t^r}{\rho^{\lambda_{ir}}} \right)^{\frac{1}{1-m_i}} \right)  \ ,
\end{equation}
being $\lambda_{ir}= N(m_i-1)+2r$.
\end{proposition}

\begin{proof}
    Assume that \eqref{condCrhoisotLr} fails. Let $\rho>0$, $t>0$ and $\sigma \in (0,1)$ be fixed and consider the cylinders $Q_\rho \subset Q_{(1+\sigma)\rho} \subset {Q}_{2\rho} \subset \Omega_T$.  Take  $\varphi= f(u) \xi^2$, for
    \[f(u) = u^{r-1} \left( \frac{(u-k)_+}{u} \right)^s \ , \quad \max \{ r-1,1 \} < s< r \ , \quad k> 0 \]
and a time-independent smooth cutoff function $\xi \in C_o^{\infty}\left( K_{(1+\sigma)\rho} \right)$ that equals one in $K_\rho$, vanishes outside $K_{(1+\sigma)\rho}$ and verifies $|D \xi| \leq {1}/{\sigma \rho} $. Arguing in the same way as in the proof of Proposition \ref{Lrbackest}, defining $\displaystyle{S_\sigma= \sup_{0 \leq \tau \leq t } \int_{K_{(1+\sigma)\rho}} u^r(x,\tau) \ dx}$, we arrive at
\begin{eqnarray*}
  \sup_{0 \leq \tau \leq t } \int_{K_\rho} u^r (x,\tau) \ dx  & \leq & \gamma(r) \int_{K_{2\rho}} u^r (x,0)  \ dx\\
  & & + \frac{\gamma(C_o,C_1,r,N)}{(r-1) \sigma^2 \rho^2} \sum_{i=1}^N  \iint_{Q_{(1+\sigma)\rho}} u^{m_i+r-1} \ \chi[u>k] \ dx d\tau \\
  & \leq & \gamma(r) \int_{K_{2\rho}} u^r (x,0)  \ dx\\
  & & + \sum_{i=1}^N  \left[ \frac{\gamma(C_o,C_1,r)}{(r-1) \sigma^2 }
    \left( \frac{t^r}{\rho^{N(m_i-1)+2r}} \right)^{\frac{1}{r}} \right] \
    \S_{\sigma}^{\frac{r-1+m_i}{r}} \\ 
     & \leq & \gamma(r)\int_{K_{2\rho}} u^r (x,0)  \ dx\\
  & & + \frac{1}{2} \S_\sigma + \sum_{i=1}^N \gamma(C_o,C_1,r, \sigma,N,m_i) \ \left( \frac{t^r}{\rho^{{N(m_i-1)+2r}}} \right)^{\frac{1}{1-m_i}} \ .
\end{eqnarray*}
By considering
\[\rho_n= \rho \sum_{j=1}^n 2^{-j} \ ; \qquad \rho_{n+1} =(1+\sigma_n) \rho_n \Rightarrow \sigma_n \geq \frac{1}{2^{n+1}}\]
and
\[ S_n = \sup_{0 \leq \tau \leq t } \int_{K_n} u^r (x,\tau) \ dx \ , \qquad \mathrm{being} \quad K_n =  \left\{ |x|<  \rho_n \right\} \  , \]
\[ S_n \leq  \frac{1}{2} S_{n+1} + b^n  \gamma(C_o,C_1,N, m_1,m_N) \left[ \int_{K_{2\rho}} u^r (x,0)  \ dx  + \sum_{i=1}^N  \left( \frac{t^r}{\rho^{\lambda_{ir}}} \right)^{\frac{1}{1-m_i}} \right]  \ , \quad b >1  \]
and now \eqref{estLrisot} is a direct consequence of Lemma \ref{iteration}.
\end{proof}
 
\begin{remark}
    The constant $\gamma$ deteriorates as $r\searrow 1$, as shown along the proof.
\end{remark}

\begin{remark}
    As a final remark, we emphasize that all the results obtained, whatever the geometric setting chosen - either the anisotropic or the standard (isotropic) one - are consistent with the theory known for a class of isotropic porous media type equations. In fact, if one considers all the exponents $m_i$, $i=1,\cdots,N$, to be equal, say 
 \[ m_i=m  \ , \qquad \forall i=1,\cdots,N \]
 the exponents appearing along the text become precisely the ones derived in the isotropic context (see for instance \cite{DGVmono} and \cite{KW}).
\end{remark}

\section*{Acknowledgements}

\noindent  S. Ciani acknowledges the support of the department of Mathematics of the University of Bologna Alma Mater and the Italian PNR (MIUR) fundings 2021-2027; E. Henriques was financed by Portuguese Funds through FCT - Funda\c c\~ao para a Ci\^encia e a Tecnologia - within the Projects UIDB/00013/2020 and UIDP/00013/2020.


\begin{thebibliography}{99}
\bibitem{AS15} S. Antontsev, S. Shmarev, {\it Evolution PDEs with nonstandard growth conditions.} Atlantis Studies in Differential Equations, 4, (2015).

       
\bibitem{AS05} S. Antontsev, S. Shmarev, {\it A model porous medium equation with variable exponent of nonlinearity: existence, uniqueness and localization properties of solutions.} Nonlinear Analysis: Theory, Methods and Applications,  60 (3), (2005), 515-545.






\bibitem{Aron85} D.G. Aronson, {\it The porous medium equation.} In: Nonlinear Diffusion Problems (Montecatini Terme, 1985), Lecture Notes in Mathematics, 1224, Springer, Berlin, 1–46, (1986).
\bibitem{Barenblatt-scaling1} G.I. Barenblatt, {\it Scaling, self-similarity, and intermediate asymptotics: dimensional analysis and intermediate asymptotics.} Cambridge University Press, 14, (1996).
\bibitem{BH80} J.G. Berryman, C.J. Holland, {\it Stability of the separable solution for fast diffusion.} Arch. Rational Mech. Anal., 74(4), (1980), 379–388. 
\bibitem{Besov} O.V. Besov, V.P. Ilin, S.M. Nikolskii, {\it Integral representations of functions and imbedding theorems.} VH Winston, Washington DC, I, 1978.


\bibitem{BF21} M. Bonforte, A. Figalli, {\it Sharp extinction rates for fast diffusion equations on generic bounded domains.} Commun. Pure Appl. Math., 74, (2021), 744–789.
\bibitem{BGV12} M. Bonforte, G. Grillo, J.L. Vazquez, {\it Behaviour near extinction for the Fast Diffusion Equation on bounded domains.} J. Math. Pures Appl., 97, (2012), 1–38.
\bibitem{CMcCS23} B. Choi, R.J. McCann, C. Seis, {\it Asymptotics Near Extinction for Nonlinear Fast Diffusion on a Bounded Domain.}
Arch. Rational Mech. Anal., 247 (16), (2023).
\bibitem{CGV22} S. Ciani, U. Guarnotta, V. Vespri,  {\it On a particular scaling for the prototype anisotropic p-Laplacian.}  In Recent Advances in Mathematical Analysis, Trends in Mathematics, Birkhauser, (2023), 289-308.

\bibitem{CSV22} S. Ciani, I. Skrypnik, V. Vespri, {\it On the Local Behavior of Local Weak Solutions to some Singular Anisotropic Elliptic Equations.} Advances in Nonlinear Analysis, 12(1), (2022), 237-265.
\bibitem{CVV} S. Ciani, V. Vespri, M. Vestberg, {\it Boundedness, Ultracontractive Bounds and Optimal Evolution of the Support for Doubly Nonlinear Anisotropic Diffusion.} Preprint (2023) (https://doi.org/10.48550/arXiv.2306.17152).

\bibitem{DB} E. DiBenedetto, {\it Degenerate Parabolic Equations}. Universitext, Springer-Verlag, New York, (1993). 
\bibitem{DGVmono} E. DiBenedetto, U. Gianazza, V. Vespri, {\it Harnack's Inequality for Degenerate and Singular Parabolic 
Equations.} Springer Monographs in Mathematics, Springer-Verlag, New York, (2012).
\bibitem{KW} E. DiBenedetto, Y.C. Kwong, {\it Intrinsic Harnack Estimates and Extinction Profile for Certain Singular Parabolic Equations.} Trans. Amer. Math. Soc., 330, (1992), 783-811.
\bibitem{DFZ20} N.M.L. Diehl, L. Fabris, J.S. Ziebell, {\it Decay Estimates for Solutions of Porous Medium Equations with Advection.} Acta Appl Math, 165, (2020), 149–162.
\bibitem{DMV} F.G. D\"uzg\"un, S. Mosconi, V. Vespri, {\it Anisotropic Sobolev embeddings and the speed of propagation for parabolic equations.} Journal of Evolution Equations, 19(3), (2019), 845-882.
\bibitem{FVV23} F. Feo, B. Volzone, J.L. Vazquez, {\it Anisotropic Fast Diffusion Equations.} Nonlinear Analysis. Theory, Methods {\&} Applications, 233, (2023), 113298.
\bibitem{FHV21} S. Fornaro, E. Henriques, V. Vespri, {\it Regularity results for a class of doubly nonlinear very singular parabolic equations.} Nonlinear Analysis, (2021), 112213.
\bibitem{EH21} E. Henriques, {\it The porous medium equation with variable exponent revisited.} J. Evol. Equ., 21 (2), (2021), 1495-1511.
\bibitem{EH11} E. Henriques, {\it Concerning the regularity of the anisotropic porous medium equation.} J. Math. Anal. Appl., 377 (2), (2011), 710-731.
\bibitem{EH08} E. Henriques, {\it Regularity for the porous medium equation with variable exponent: The singular case.} J. Differ. Equations, 244 (10), (2008), 2578-2601.
\bibitem{HU06} E. Henriques, J.M. Urbano, {\it Intrinsic scaling for PDE’s with an exponential nonlinearity.} Indiana Univ. Math. J., 55 (5), (2006), 1701-1721.
\bibitem{KL08} J. Kinnunen, P. Lindqvist, {\it Definition and properties of supersolutions to the porous medium equation.} J. Reine Angew. Math., 618,  (2008), 135–168.
\bibitem{LWW07} W. Liu, M. Wang, B. Wu, {\it Extinction and decay estimates of solutions for a class of porous medium equations.} Journal of Inequalities and Applications 2007, (2007), 1-8.
\bibitem{MNS23} M. Misawa, K. Nakamura, M.A.H Sarkar, {\it A finite time extinction profile and optimal decay for a fast diffusive doubly nonlinear equation.} Nonlinear Differ. Equ. Appl., (2023), 30-43.
\bibitem{Porzio23} M. M. Porzio, {\it On the speed of decay of solutions to some partial differential equations.} J. Math. Anal. Appl. (2023), 127535.
\bibitem{Sabinina62} E.S. Sabinina, {\it On a class of non-linear degenerate parabolic equations.} Dokl. Akad. Nauk SSSR, 143, (1962), 794-797.
\bibitem{SongJian05} B. H. Song, H. Y. Jian, {\it Fundamental solution of the anisotropic porous medium equation.} Acta Math. Sinica, 21 (5), (2005), 1183–1190.
\bibitem{SongJian06} B. H. Song, H. Y. Jian, {\it Solutions of the anisotropic porous medium equation in $\R^N$ under an $L^1$ initial value.} Nonlinear Analysis, 64, (2006), 2098 – 2111.
\bibitem{Songexist01} B. Song, {\it  Anisotropic diffusions with singular advections and absorptions, Part I, Existence.} Appl. Math. Lett., 14, (2001), 811–816.
\bibitem{Songuniq01} B. Song, {\it  Anisotropic diffusions with singular advections and absorptions, Part II, Uniqueness.} Appl. Math. Lett., 14, (2001), 817–823.
 \bibitem{Troisi} M. Troisi, {\it Teoremi di inclusione per spazi di Sobolev non isotropi.} Ricerche Mat 18(3), (1969), 3-24.
 \bibitem{Vaz12} J.L. Vazquez, {\it The porous medium equation: mathematical theory.} Oxford Mathematical Monographs, Oxford Science Publications, Clarendon Press, Oxford, (2012).
\bibitem{Vaz06} J.L. Vazquez, {\it Smoothing and decay estimates for nonlinear diffusion equations of porous medium type.} Oxford Lecture Series in Mathematics and its Applications, Oxford University Press, 33, (2006)
\bibitem{ZK50} Y.B. Zel‘dovič, A.S. Kompaneec, {\it On the theory of propagation of heat with the heat conductivity depending upon the temperature.} In: Collection in Honor of the Seventieth Birthday of Academician A. F. Ioffe, Izdat. Akad. Nauk SSSR, Moscow, (1950), 61–71.

\end{thebibliography}
\end{document}